\pdfoutput=1
\RequirePackage{ifpdf}
\ifpdf 
\documentclass[pdftex]{sigma}
\else
\documentclass{sigma}
\fi

\numberwithin{equation}{section}

\newtheorem{Theorem}{Theorem}[section]
\newtheorem{Corollary}[Theorem]{Corollary}
\newtheorem{Lemma}[Theorem]{Lemma}
\newtheorem{Proposition}[Theorem]{Proposition}
 { \theoremstyle{definition}
\newtheorem{Definition}[Theorem]{Definition}

 }

\usepackage{tabu}

\renewcommand{\epsilon}{\varepsilon}
\newcommand{\vep}{\varepsilon}
\newcommand{\wt}{\widetilde}

\newcommand{\del}{\partial}
\newcommand{\RenV}{\mathrm{RenV}}
\newcommand{\olg}{\overline{g}}
\newcommand{\olG}{\overline{G}}
\newcommand{\olGamma}{\overline{\Gamma}}

\newcommand{\olgin}{\overline{\gin}}

\newcommand{\olE}{\overline{E}}
\newcommand{\frakc}{\mathfrak c}
\newcommand{\calO}{\mathcal O}
\newcommand{\calA}{\mathcal A}
\newcommand{\calU}{\mathcal U}
\newcommand{\calC}{\mathcal C}
\newcommand{\e}{\epsilon}
\newcommand{\CC}{\mathbb C}

\newcommand{\calL}{\mathcal L}
\newcommand{\gin}{{g^0}}

\newcommand{\tr}{\mathrm{tr}}

\newcommand{\gbar}{\overline{g}}

\newcommand{\Mbar}{\overline{M}}

\renewcommand{\hbar}{\overline{h}}

\newcommand{\FP}{\underset{z=0}{\operatorname{FP}}\,}
\newcommand{\Res}{\underset{z=0}{\operatorname{Res}}\,}
\DeclareMathOperator{\Rc}{Rc}

\newcommand{\Renint}{\mathop{\vphantom{T}}^R \! \! \! \int}

\DeclareMathOperator{\hess}{Hess}
\DeclareMathOperator{\re}{Re}

\begin{document}
\allowdisplaybreaks

\newcommand{\arXivNumber}{1607.08558}

\renewcommand{\PaperNumber}{057}

\FirstPageHeading

\ShortArticleName{Ricci Flow and Volume Renormalizability}

\ArticleName{Ricci Flow and Volume Renormalizability}

\Author{Eric BAHUAUD~$^\dag$, Rafe MAZZEO~$^\ddag$ and Eric WOOLGAR~$^\S$}

\AuthorNameForHeading{E.~Bahuaud, R.~Mazzeo and E.~Woolgar}

\Address{$^\dag$~Department of Mathematics, Seattle University, 901 12th Ave, Seattle, WA 98122, USA}
\EmailD{\href{mailto:bahuaude@seattleu.edu}{bahuaude@seattleu.edu}}

\Address{$^\ddag$~Department of Mathematics, Stanford University, Stanford, CA 94305, USA}
\EmailD{\href{mailto:mazzeo@math.stanford.edu}{mazzeo@math.stanford.edu}}

\Address{$^\S$~Department of Mathematical and Statistical Sciences, University of Alberta,\\
\hphantom{$^\S$}~Edmonton, Alberta, T6G 2G1, Canada}
\EmailD{\href{mailto:ewoolgar@ualberta.ca}{ewoolgar@ualberta.ca}}

\ArticleDates{Received December 06, 2018, in final form July 30, 2019; Published online August 07, 2019}

\Abstract{With respect to any special boundary defining function, a conformally compact asymptotically hyperbolic metric has an asymptotic expansion near its conformal infinity. If this expansion is even to a certain order and satisfies one extra condition, then it is possible to define its renormalized volume and show that it is independent of choices that preserve this evenness structure. We prove that such expansions are preserved under normalized Ricci flow. We also study the variation of curvature functionals in this setting, and as one application, obtain the variation formula
\[
\frac{{\rm d}}{{\rm d}t} {\rm RenV}\big(M^n, g(t)\big) = -\mathop{\vphantom{T}}^R \! \! \! \int_{M^n} (S(g(t))+n(n-1) ) {\rm d}V_{g(t)},
\]
where $S(g(t))$ is the scalar curvature for the evolving metric $g(t)$, and $\mathop{\vphantom{T}}^R \! \! \! \int (\cdot) {\rm d}V_g$ is Riesz renormalization. This extends our earlier work to a broader class of metrics.}


\Classification{53C44}

\section{Introduction}

A basic invariant of a compact Riemannian manifold is its volume. Remarkably, certain conformally compact
asymptotically hyperbolic metrics, including Poincar\'e--Einstein metrics in even dimensions, have a well-defined,
finite renormalized volume. In this paper we study the behavior of this quantity under normalized Ricci flow.

Conformally compact asymptotically hyperbolic (AH) spaces $(M^n, g)$ are a class of complete Riemannian manifolds modeled on the Poincar\'e disk model of hyperbolic space. Any such metric is defined and complete on the interior of a~compact manifold with boundary $\Mbar = M \cup \partial M$, and takes the form $g = \rho^{-2} \gbar$ where $\gbar$ is some smooth metric on $\Mbar$ and where $\rho$ is a smooth boundary defining function. We require that $|{\rm d}\rho/\rho|^2_g \to 1$ at $\partial M$, which ensures that the sectional curvatures of $g$ tend to $-1$ near $\partial M$. The conformal class $[\rho^2 g |_{T\partial M}]$ is a well-defined conformal class on~$\partial M$, and $\del M$ with this conformal class is called the conformal infinity of~$(M,g)$. This class of metrics provides the setting for many interesting problems in geometric analysis and physics.

This class contains the `special' Poincar\'e--Einstein (PE) metrics, which by definition are the conformally compact metrics which satisfy
\begin{gather} \label{e1}
E := \Rc(g) + (n-1)g = 0,\qquad n = \dim M.
\end{gather}
There are differences in various parts of this theory depending on the parity of $n$, but in this paper we restrict ourselves
{\it entirely} to the case where $n$ is even. We emphasize our convention that $n = \dim M$ is the dimension of the bulk manifold,
which is in contrast to certain other papers on PE metrics where $n$ denotes the dimension of the conformal boundary instead.
It turns out that a choice of representative of the conformal infinity uniquely determines a boundary defining function and
a diffeomorphism of a collar neighbourhood of $\partial M$ with $\partial M \times (0,1)$. In terms of this identification,
there is a formal series solution $g$ to equation~\eqref{e1} called the Fefferman--Graham expansion. This expansion is
even up to order $n$ with respect to the special defining function, with coefficients locally determined by the choice of metric
in the conformal infinity. Furthermore, the first odd term in the expansion has vanishing trace. Using these properties,
one may then compute the volumes of compact truncations in the exhaustion of the manifold determined by these special
boundary defining functions. Discarding the divergent terms in the expansion of this one-parameter family of volumes
defines the renormalized volume. Considerations of parity in this expansion then shows that this renormalized volume
is independent of the choice of representative of the conformal infinity.

We study here this renormalized volume in a slightly more general setting, where the asymptotic expansion of the metric
has qualitatively similar features. The first possibility is to relax~\eqref{e1} and simply require that $g$ is asymptotically
Poincar\'e--Einstein (APE) in the sense that $|E|_g = O(\rho^n)$ for some (hence any) boundary defining function $\rho$. In
this case one has the same Fefferman--Graham expansion of the metric and the same volume renormalization scheme as
in the PE case. However, we relax these conditions even further, and require only that the expansion of the metric be even
to a critical order and that the first odd term have vanishing trace. We call AH metrics satisfying the first condition {\it
partially even}, and metrics satisfying both conditions {\it volume renormalizable}. Thus both PE metrics and APE metrics are
volume renormalizable. However, no special properties are required of the nonzero coefficients in the expansion of
a volume renormalizable metric, and hence the $E$ tensor need only vanish to second order.

AH metrics with various assumptions on the evenness of the expansion have been considered frequently before,
see \cite{Guillarmou, MM} and the more recent work by Vasy \cite{V} regarding the role of evenness in establishing
the meromorphic extension and properties of the resolvent of the Laplacian of an AH metric.

We also consider the normalized Ricci flow
\begin{gather} \label{normRF}
\partial_t g = -2 ( \Rc(g) + (n-1)g ), \qquad t \in [0,T), \qquad g(0) = \gin.
\end{gather}
The $(n-1) g$ term ensures that the conformal infinity is fixed in time. Observe that hyperbolic metrics and more generally
PE metrics are stationary points of this flow.

The preservation of the APE condition under the normalized Ricci flow and the variation of renormalized volume was studied
in our earlier paper~\cite{BMW}. Here we turn to these questions in the more general class of partially even and volume renormalizable
metrics. We now state our main results.

\begin{Theorem}\label{TheoremA} Suppose that $\big(M,\gin\big)$, $\dim M =: n = 2m$, is partially even, and let $g(t)$ be a~solution of~\eqref{normRF} with $g(0) = \gin$
with maximal interval of existence $[0,T_0)$. Then $(M, g(t))$ remains partially even for $t < T_0$. If, furthermore, $\big(M,\gin\big)$ is volume renormalizable then $(M, g(t))$ remains volume renormalizable for $t < T_0$.
\end{Theorem}

Our second main result computes the variation of the renormalized volume along this flow. We cast this slightly more generally
by proving that the variation through volume renormalizable metrics of any curvature functional is just the renormalized first
variation of the same form that, on compact manifolds, yields the Euler--Lagrange equations, see Theorem~\ref{thm:renorm-curv-fcl} below.
An immediate application is the
\begin{Theorem} \label{TheoremB}
Let $\big(M,\gin\big)$, $n = 2m$, be volume renormalizable, and $g(t)$ the solution of~\eqref{normRF} with $g(0) = \gin$. Then
\begin{gather}\label{eq0.3}
\frac{{\rm d}}{{\rm d}t} \RenV = - \Renint_M S(g(t)) + n(n-1) \, {\rm d}V_{g(t)} .
\end{gather}
\end{Theorem}
Here $\Renint (\cdot) {\rm d}V_{g}$ indicates that the integral is not classically convergent, but must be regularized in a way to be
made precise below. The analogue of formula~(\ref{eq0.3}) in~\cite{BMW} appears almost exactly the same except that when
$g(t)$ is APE, the integral is convergent and there is no need to renormalize the integral.

Finally let us present a nonexhaustive list of related works studying the Ricci flow in the context of AH metrics.
Short-time existence was obtained in \cite{WXShi} and a uniqueness result in~\cite{Chen-Zhu}. Preservation of
conformally compact metrics along the flow was obtained in \cite{Bahuaud}. Stability of the flow around hyperbolic space
or other PE metrics under various hypotheses have been obtained in \cite{Bamler,QSW,SSS}. A long-time existence and
convergence result for rotationally symmetric AH metrics was obtained in \cite{BahuaudWoolgar}. The present authors
studied the behavior of APE metrics along the flow in~\cite{BMW}. The evolution of the mass aspect tensor was obtained in~\cite{BW}. In \cite{R} it is proved that polyhomogeneity is preserved along a flow adapted to asymptotically complex hyperbolic metrics. More generally still, \cite{Ammar} proves that the polyhomogeneity of a metric with a Lie structure fibred at infinity is preserved along the Ricci flow.

This paper is organized as follows. In Section~\ref{FGofAH} we review the families of asymptotically hyperbolic metrics
considered here and establish some of basic properties involving even expansions. In Section \ref{PEandRF} we prove that
the normalized Ricci flow of a partially even metric remains partially even. In Section \ref{VRandRF} we complete the proof
of Theorem~\ref{TheoremA} by showing that the volume renormalizability condition also persists under the Ricci flow. We
defer several long computations needed in this section to Appendix~\ref{appendix}. Finally in Section~\ref{VarofCurvFcl} we review the
Riesz renormalization and apply it to the variation of curvature functionals. We then prove Theorem~\ref{TheoremB}.

\section{Fefferman--Graham expansions of AH metrics} \label {FGofAH}

\subsection{Taxonomy of conformally compactifiable metrics}
As in the introduction, let $\Mbar$ be a smooth compact $n$-manifold with boundary and $M$ its interior. We assume $n$ is even.

A metric $g$ on $M$ is called conformally compact if
\begin{gather*}
g = \rho^{-2} \olg ,
\end{gather*}
where both $\rho$ and $\olg$ are $\calC^\infty$ up to $\del M$, and $\olg$ is a metric on $\Mbar$.

Any such metric has a sort of boundary value on $\del M$,
\begin{gather*}
\mathfrak c(g) = \big[ \rho^2 g \big|_{T\del M} \big],
\end{gather*}
called its conformal infinity. Furthermore, a conformally compact metric is called \emph{asymptotically hyperbolic},
abbreviated as AH, if $|{\rm d}\rho|^2_{\olg} = 1$ at $\del M$. Note that $g$ is unchanged if we replace $\rho$ and $\olg$
by $a\rho$ and $a^2 \olg$ for any smooth positive function $a$, but nonetheless both $\mathfrak c(g)$ and the AH condition are well-defined.

Now, if $g$ is AH and we select a representative metric $h_0 \in \frakc(g)$, then there exists a uniquely determined
`special' boundary defining function $x$ which satisfies $|{\rm d}x|^2_{x^2 g} \equiv 1$ in a~neighborhood~$V_x$ in $\Mbar$ of the boundary and $x^2 g|_{T \del M} = h_0$.
With this choice of $x$, the flow lines of the gradient $\nabla^{\olg} x$ allow us to identify $V_x$ with a neighborhood of $\del M$ of the form $U_x = [0,\vep_x) \times \del M$, via a diffeomorphism
\[ \phi\colon \ U_x \longrightarrow \phi(U_x) \subset V_x, \]
where $\phi(s,y)$ is the flow of $\nabla^{\olg} x$ for $s$ units of time. By Gauss's lemma
\begin{gather} \label{eq1.3}
\phi^*g = \frac{{\rm d}s^2 + h(s)}{s^2},
\end{gather}
where $h(s)$ is a smooth family of metrics on $\del M$. This is the Graham--Lee normal form~\cite{GL}, and we follow common practice of writing~$g$ and $x$ in place of $\phi^* g$ and $s$ above. We shall be primarily concerned with the terms in the asymptotic expansion
\begin{gather*}
h(x) = h_0 + h_1 x + h_2 x^2 + \cdots.
\end{gather*}
For some of the computations below we will also need a full coordinate system near the boundary of~$M$. Let $\{y^{\alpha}\}$ be any coordinates on~$\partial M$ extended to be constant along the integral curves of~$\nabla^{\olg} x$. We use Greek indices to index these tangential directions, $x$ to index $\frac{\partial}{\partial x}$, and Latin indices to include both tangential and normal directions. In the expansion above, the coefficient tensors are functions of the $y$-coordinates alone.

\begin{Definition}\label{definition1.1} \label{def1} An AH metric $g$ is called \emph{even to order $2\ell$} if, in Graham--Lee normal form~\eqref{eq1.3}, the expansion for~$h(x)$ contains no terms $h_{2j+1}$ with $j \leq \ell-1$, i.e.,
\begin{gather*}
h(x) \sim h_0+x^2h_2 +\dots +x^{2\ell}h_{2\ell} +x^{2\ell+1}h_{2\ell+1}+ \cdots.
\end{gather*}
Our main focus will be on AH metrics $g$ which are even to order $n-2$, and we call these simply \emph{partially even}.
\end{Definition}

That this notion and the next are well defined, independent of the choice of special defining function is proved in Proposition~\ref{proposition1.5} below.

\begin{Definition}\label{definition1.2}We say that $g$ is \emph{volume renormalizable} (VR) if it is partially even and in addition satisfies
\begin{gather*}
\tr^{h_0} h_{n-1} = 0.
\end{gather*}
\end{Definition}

For volume renormalizable metrics, one can compute the volume in the compact set $\{x \geq \e\}$, expanding the result in powers of $\epsilon$ (there can also be $\log \epsilon$ terms when the bulk dimen\-sion~$n$ is odd) and taking the limit at $\e = 0$ after discarding the singular terms in this expansion. While this may be done for any AH metric with a fixed special defining function $x$, this has an invariant meaning which is independent of the choice of $x$ precisely when the metric is VR~\cite{Gr}.

We now recall one of the most important classes of AH metrics
\begin{Definition}\label{definition1.3} An AH metric $g$ is called \emph{Poincar\'e--Einstein} (PE) if its modified Ricci tensor~$E(g)$, defined in~\eqref{e1}, vanishes identically.
\end{Definition}
Since our primary concern is with asymptotic expansions at the boundary, it is also natural to introduce the
\begin{Definition}\label{definition1.4}
An AH metric $g$ is called \emph{asymptotically Poincar\'e--Einstein} (APE) if $|E(g)|_{g} = \calO( x^n)$.
\end{Definition}

If $g$ is APE, or in particular PE, then it is known, see~\cite{FG2}, that every odd coefficient $h_{2j+1}$, $0 \leq j \leq (n-4)/2$, must
vanish while the even coefficients $h_{2j}$, $1 \leq j \leq (n-2)/2$, can be expressed as universal differential operators applied
to $h_0$, see \cite{FG2}, and finally, $\tr^{h_0} h_{n-1} = 0$. This implies that when $n$ is even,
\begin{gather*}
	\text{PE $\implies$ APE $\implies$ VR $\implies$ partially even $\implies$ AH.}
\end{gather*}

We have not included the slightly different results when $n$ is odd; the key difference in that case is that when $g$ is PE, the expansion for $h(x)$ may also include the terms $x^{n-1+j} (\log x)^\ell \tilde{h}_{j\ell}$, $j \geq 0$. When $n=2m$, if $g$ is merely $\calC^2$ conformally compact and Einstein, then a regularity theorem due to Chru\'sciel, Delay, Lee and Skinner~\cite{CDLS} shows that if $h_0 \in \calC^\infty$, then $h(x)$ is automatically smooth up to $x=0$, i.e., it is not necessary to assume a priori that $g$ is smoothly conformally compact, so long as its conformal infinity contains a smooth representative.

\subsection{Partially even metrics}

We now discuss the extent to which even expansions are well-defined. The procedure leading to Graham--Lee normal form discussed in the previous section establishes a bijection between representatives of the conformal infinity $\mathfrak c(g)$ and the set of special defining functions. The following lemma is due to Graham \cite{Gr} with an important addition by Guillarmou \cite{Guillarmou}.

\begin{Proposition}\label{proposition1.5}
If an AH metric $g$ is even to order $2\ell$ in the normal form corresponding to one choice of metric $h_0 \in \frakc (g)$,
then it is even to the same order with respect to any other metric $h_0' \in \frakc(g)$. Further, if $x$ and $x'$ are the special boundary defining functions corresponding to $h_0$ and $h_0'$ respectively, and $y = y^{\alpha}$ and $y' = {y'}^{\alpha}$ are corresponding choices of smooth coordinates on $\partial M$ extended to be constant along the respective gradient flow lines, then coordinate change of $(\phi')^{-1}\circ \phi$ on $U_x$ has expansions
\[ x' = x \sum_{j=0}^{\ell+1} a_j(y) x^{2j} + O\big(x^{2\ell+4}\big), \qquad y' = \sum_{j=0}^{\ell+1} b_j(y) x^{2j} + O\big(x^{2\ell+3}\big), \]
for smooth functions $a_j$ and $b_j$ on $\partial M$.
\end{Proposition}
\begin{proof}[Sketch of proof] The first part of this lemma is proved in Graham \cite{Gr} by first writing $h_0' = {\rm e}^{2\omega_0} h_0$, where $\omega_0 \in \calC^\infty(\del M)$. The special boundary defining function $x'$ associated to $h_0'$ is then given by $x' = {\rm e}^{\omega} x$, where $\omega$ satisfies $\big| {\rm d}({\rm e}^\omega x)/ ({\rm e}^\omega x)\big|^2_g \equiv 1$, $\omega|_{x=0} = \omega_0$. Expanding the eikonal equation, we obtain
\begin{gather}\label{eq1.7}
2x \langle {\rm d}\omega, {\rm d}x \rangle_{\olg} + x^2 |{\rm d}\omega|^2_{\olg} = 0, \qquad \omega |_{\del M} = \omega_0.
\end{gather}
Inductively computing derivatives of $\omega$ and using the parity of $g$ allows one to prove that~$\omega$ is even in~$x$ to order~$2 \ell+2$.

The second part of this lemma is proved by Guillarmou~\cite{Guillarmou} by explicitly writing out the gradient flow equations of $x'(x,y)$ and $y'(x,y)$ and again arguing by induction using the parity of $\omega$ and~$h(x)$.
\end{proof}

It follows from this proof that if $g$ is an AH metric which is even to order $2\ell$, then it defines an equivalence class of defining functions $[x]$, where $x \sim x'$ if $x'/x$ is even to order $2\ell+2$. This motivates another approach to studying even expansions: we may translate parity conditions for the metric into regularity statements about extensions over the boundary. To this end, let~$\wt{M}$ denote the double of $M$ across its boundary. The
atlas of $\calC^\infty$ functions on $\wt{M}$ is determined uniquely once we specify an identification of a collar neighborhood
$\calU$ of $\del M$ with the product $\del M \times [0,\epsilon)$; the double of this neighborhood, $\wt{\calU}$, is then
endowed with the atlas of $\calC^\infty$ functions on the product $\del M \times (-\epsilon, \epsilon)$. While this atlas a priori depends on the choice of $x$, using Proposition~\ref{proposition1.5} we thus see for an AH metric even to order $2 \ell$ that there is a well defined $\calC^{2\ell+2}$ structure on the doubled manifold $\wt{M}$.

Having checked that the notion of even expansion is well defined, we now check that the volume renormalizability condition is well defined.

\begin{Proposition}\label{proposition1.5b}
If an AH metric $g$ is a volume renormalizable metric in the normal form corresponding to one choice of metric $h_0 \in \frakc (g)$,
then it is volume renormalizable to the same order with respect to any other metric $h_0' \in \frakc(g)$.
\end{Proposition}
\begin{proof} Choose $h_0 \in \frakc (g)$ and obtain a special defining function $x$ so that the metric $g$ may be written in the normal form of equation~\eqref{eq1.3}, and even to order $n-2$, with $\tr^{h_0} h_{n-1} = 0$. For any $h_0' \in \frakc(g)$, Proposition~\ref{proposition1.5} guarantees that the expansion of $g$ relative to $x' = {\rm e}^{\omega} x$ remains even to order $n-2$. We need only check the vanishing trace condition, $\tr^{h_0'} h_{n-1}' = 0$. Since $h_0$ and $h_0'$ are conformal, note that $\tr^{h_0'} h_{n-1} = 0$.
	
Consider the auxiliary product metric $G' = ({\rm d}x')^2 + h_0'$ on $U_{x'}$, and define a function $t' = \tr^{G'} \big[(x')^2 g\big]$. By construction the $(n-1)$-st term in the expansion of $t'$ relative to $x'$ is $\tr^{h_0'} h_{n-1}'$. On the other hand, the function $t'$ may be alternatively computed by
\begin{gather*}
	t' = \tr^{G'} \big[(x')^2 g\big] = \tr^{G'} \left[ \frac{(x')^2}{x^2} \gbar \right]
	= {\rm e}^{2\omega} \tr^{G'} [ \gbar ] = {\rm e}^{2\omega} \tr^{G'} \big[ {\rm d}x^2 + h_{\alpha \beta}(x,y) {\rm d}y^{\alpha} {\rm d}y^{\beta} \big].
\end{gather*}
Now express $G'$ in $(x,y)$ coordinates using the expansions for $x'$ and $y'$ in terms of $x$ of Proposition~\ref{proposition1.5}, and combine with the even expansion for $h(x)$, the fact $\gbar$ contains no cross terms and the fact that $\tr^{h_0'} h_{n-1} = 0$ to obtain that $\tr^{G'} [ \gbar ]$ is even in~$x$ to order~$n$. Using the parity for $\omega$ we then find $t'$ is even to order $n$ in $x$ and thus defines a $\calC^n$ function on the doubled manfiold. Thus the $(n-1)$-st term in the expansion of $t'$ relative to~$x'$ vanishes, as required.
\end{proof}

We now extend this discussion by considering how to recognize partially even AH metrics which are not written in
normal form. Thus fix an AH metric $g$ and suppose that $x$ is an arbitrary boundary defining function. As before
we set $\olg = x^2 g$ and use the gradient lines of $\nabla^{\olg} x$ to define a~(pointed) smooth structure $\calC^{\infty}_x$ on the
double~$\wt{M}$, and finally, consider the extended metric~$\tilde{g}$ which restricts to $\olg$ on $M$ and
satisfies $I_x^* \tilde{g} = \tilde{g}$, where $I_x\colon \wt{M} \to \wt{M}$ is reflection across the submanifold $\del M \subset \wt{M}$.
\begin{Definition}\label{definition1.6} \label{defn-even}
We say that $g$ is even to order $j$ relative to an arbitrary boundary defining function $x$ if $\tilde{g}$ is $\calC^{j}_x$ on~$\wt{M}$ (still recalling that it is $\calC^\infty$ on each side of $\del M$).
\end{Definition}

\noindent To explain the definition further, suppose $\{y^{\alpha}\}$ is any choice of coordinates on $\partial M$ extended as usual to be constant along the integral curves of $\nabla^{\gbar} x$. Consider a symmetric $2$-tensor $T$
\[ T= T_{xx}(x,y) {\rm d}x^2+ T_{x\alpha}(x,y) {\rm d}x {\rm d}y^{\alpha} + T_{\alpha \beta}(x,y) {\rm d}y^{\alpha} {\rm d}y^{\beta}, \]
defined and smooth for $x \geq 0$. The $I_x$ invariant extension is then given by the same formula for $x \geq 0$ and for $x < 0$,
\[ \widetilde{T} = T_{xx}(-x,y) {\rm d}x^2 - T_{x\alpha}(-x,y) {\rm d}x {\rm d}y^{\alpha} + T_{\alpha \beta}(-x,y) {\rm d}y^{\alpha} {\rm d}y^{\beta}. \]

Now suppose $j = 2\ell$ and that $\widetilde{T}$ is $\calC^{j}_x$. Then the component functions $T_{xx}$ and $T_{\alpha \beta}$ are even to order $2\ell$ while $T_{x \alpha}$ is odd to order~$2\ell+1$ (i.e., it is of the form $x m_{x\alpha}$ where $m_{x\alpha}$ is even to order~$2\ell$). To explain this last condition,
note that a term of order $2\ell$ in the expansion for~$T_{x\alpha}$ with the correct parity must be of the form
$a(y) (\mathrm{sgn}\, x) |x|^{2\ell}$, and hence is not $\calC^{2\ell}_x$. On the other hand, if $j = 2\ell+1$, then evenness to
order $j$ implies that $T_{xx}$ and $T_{\alpha \beta}$ are even to order $j+1=2\ell+2$ while $T_{x\alpha}$ is still just
odd to order $2\ell+1$. Moreover, if $T$ is a metric and the components of $T$ satisfy these parity conditions, then by Cramer's formula, the
components of the inverse metric are also $\calC^{2\ell}_x$ across $x=0$, hence these too satisfy the same parity conditions.

Returning to metrics, clearly, if $g$ is even to order $2\ell$ in the sense of Definition~\ref{definition1.1}, then $g$ is even to order $2\ell$ in the sense of Definition~\ref{definition1.6} relative to a special defining function~$x$. As the following result shows, the converse is true too.

\begin{Lemma}\label{lemma1.7} The metric $g$ is even to order $2\ell$ relative to any fixed boundary defining function $($cf.\ Definition {\rm \ref{definition1.6})} if and only if $g$ is even to order $2\ell$ with respect to any special boundary defining function $($cf.\ Definition~{\rm \ref{definition1.1})}.
\end{Lemma}

\begin{proof}If $x$ denotes the initial defining function, let $h_0 = x^2 g$ restricted to $T\del M$. We seek $x' = {\rm e}^{\omega} x$ such that $ |{\rm d}x'/x'|^2_g = 1$, $\omega_0(y) = \omega(0,y) \equiv 0$, or equivalently
\begin{gather*} 
2 \langle {\rm d}x, {\rm d}\omega\rangle_{\olg} + x |{\rm d}\omega|^2_{\olg} = \frac{1 - |{\rm d}x|^2_{\olg}}{x}, \qquad \omega |_{\del M} = 0.
\end{gather*}
In contrast to equation~\eqref{eq1.7}, the right hand side of this equation is nonvanishing, but there still exists a unique solution $\omega$ which is $\calC^\infty$ for $0 \leq x < \epsilon$, and the only issue is the degree of smoothness of its even extension across $x=0$. However, since $g$ is even and asymptotically hyperbolic, $ x^{-1}\big(|{\rm d}x|^2_{\olg} - 1\big) = x^{-1} (\olg^{xx} - 1)$ is odd to order $2\ell-1$, hence lies in $\calC^{2\ell-1}_x$, and the solution of this equation is one order smoother, i.e., lies in $\calC^{2\ell}_x$. This proves that $\omega$ is vanishes to second order and is even in $x$ to order~$2\ell$, and thus $x'$ has an odd expansion in~$x$ to order~$2 \ell+1$.

Now set $g' = (x')^2 g = {\rm e}^{2\omega} \olg$ and write $\tilde{g}'$ for its $I_x$-invariant extension. From the paragraph above, $\tilde{g}' \in \calC^{2\ell}_x$. In fact, the metric is $\calC^\infty$ in the tangential direction and its irregularity is only in the direction normal to the boundary (technically, it is polyhomogeneous at the boundary). It is standard to show that the exponential mapping $\Phi\colon NM \to M$ from the normal bundle of~$\del M$ to~$M$, $\Phi(x',y) = \exp_y( x' \nu(y))$ is smooth in the `base' variable $y \in \del M$ and $\calC^{2\ell+1}_x$ in~$x$, by writing out the flow equations as a first-order system in a way similar to \cite[p.~6]{Guillarmou}. Hence $\Phi^*\tilde{g}' = {\rm d}(x')^2 + \tilde{h}'_{\alpha \beta}(x',y) {\rm d}y^\alpha {\rm d}y^\beta$, where $\tilde{h}'_{\alpha \beta} \in \calC^{2\ell}_x$, or finally, $g'$ is even to order $2\ell$ in~$x$, and thus in~$x'$.
\end{proof}

As a consequence of this lemma, an AH metric that is even to order $2\ell$ in the sense of Definition~\ref{definition1.6} gives rise to a well-defined $\calC^{2\ell}$-structure on the doubled manifold. This structure contains the $C^{2\ell+2}$-structure defined by special defining functions discussed in the previous section.

\section{Partially even metrics and Ricci flow} \label{PEandRF}

Our aim in this section is to prove the following.

\begin{Proposition} \label{prop:even-expansions-preserved}	Suppose that $\gin$ is an AH metric which is even to order $2\ell$. Let $g(t)$ be the unique solution to the Ricci flow	equation \eqref{normRF} with initial condition $\gin$. Then $g(t)$ remains even to order $2\ell$ throughout its time of existence.
\end{Proposition}

Notice that we do not use the gauged Ricci flow equation here since we are taking the existence of the solution as given. Indeed, we take from \cite{Bahuaud} the fact that there is a unique solution to~\eqref{normRF} and this remains smoothly conformally compact throughout its maximal time interval of existence.

The proof of the proposition requires some further preliminaries. First, we choose a representative~$h_0$ of the conformal infinity of the initial metric $\gin$. Obtain a special defining function~$x$ corresponding to this choice. We now fix this~$x$ and use it for all our computations. Note that as the metric~$g(t)$ evolves under the Ricci flow $x$ will no longer be a special defining function, so we will need to consider the extended notion of even metric from the previous section.

To proceed we make some remarks concerning an AH metric $g = g(t)$ at fixed time. Appealing to the transformation formula for the Ricci tensor under a conformal change of metric, see~\cite{Besse}:
if $g = x^{-2} \olg$, then $E(g) = \Rc(g) + (n-1)g = x^{-2} \olE(\olg)$, where
\begin{gather}\label{eq2.1}
\olE(\olg) = - (n-1) \big(|{\rm d}x|^2_{\gbar} -1\big) \olg + ( (n-2) \hess_{\olg} (x) + (\Delta_{\olg} x) \olg ) x + \Rc(\gbar) x^2 .
\end{gather}

We shall analyze the equation
\begin{gather}
\del_t \olg = -2 \olE( \olg),\label{rf1}
\end{gather}
which is obtained if we multiply \eqref{normRF} by $x^2$.

As a first step, consider the restriction of this equation to $x = 0$:
\begin{gather*}
 \del_t (\olg |_{x=0} ) = 2(n-1) \big( |{\rm d}x|_{\olg}^2 - 1\big) (\olg |_{x=0} ) .
\end{gather*}
Since $g$ is AH, the leading coefficient $|{\rm d}x|^2_{\olg} - 1$ vanishes at $t=0$, which shows that the restriction (not just the
pullback) of $\olg$ to the boundary is invariant under the flow. It is thus reasonable to write the evolving solution
in the form $g(t) = \gin + k(t)$, where $|k(t)|_{\gin} = \calO(x)$. This was the ansatz in~\eqref{normRF} in~\cite{Bahuaud}, and
a uniqueness theorem in this setting~\cite{Chen-Zhu} shows that any solution of~\eqref{rf1} in this quasi-isometry class
must agree with the solution in~\cite{Bahuaud}. The proof in~\cite{Bahuaud} also shows that $k(t) \in \calC^\infty(\overline{M})$
for every $t \geq 0$. Consequently, $\olE\big( x^2 g(t)\big) = 0 $ at $\del M$ for all $t \geq 0$.

Now rewrite \eqref{rf1} as an equation for the evolving tensor $v(t) = \olg(t) - \olgin$:
\begin{gather}
 \del_t v(t) = 2 (n-1) \big(|{\rm d}x|^2_{\gbar(t)} -1\big) \olg(t) -2 ( (n-2) \hess_{\olg(t)} (x) + (\Delta_{\olg(t)} x) \olg(t) ) x\nonumber\\
\hphantom{\del_t v(t) =}{} -2 \Rc(\gbar(t)) x^2, \label{rf2}
\end{gather}
where here $v \in x C^{\infty}(\Mbar)$, which we expand as $v \sim \sum \{ v(t,y) \}_n x^n$.

We introduce the following notation for coefficients in smooth power series expansions. Suppose that $f(x,y) \sim a_0(y) + a_1(y) x + \cdots$ is smooth. To simplify notation below, we denote by~$\{f\}_{n}$ the $n^{\mathrm{th}}$ coefficient function $a_n(y)$.

\begin{proof}[Proof of Proposition~\ref{prop:even-expansions-preserved}]
We now commence an inductive proof to show that if $\olgin$ is even to order $2\ell$ then $\olg(t)$ remains even to order $2\ell$. We assume $\ell \geq 1$. Note that we have relegated most of the long parity computations to the appendix.

\subsection*{Base case}
To check the base case, we must check $\olg(t)$ is even to order $2$ given that $\olgin$ is even to order $2$. This means we must check that along the flow both
\begin{align*}
v_{xx}, \ v_{\alpha \beta}\ \mbox{remains even to order} \ 2, \ \mbox{and} \ v_{x\alpha} \ \mbox{remains odd to order} \ 3.
\end{align*}
To do this we check the coefficients $\{v_{xx}\}_1$, $\{v_{\alpha \beta}\}_1$ and $\{v_{x\alpha}\}_2$ all vanish along the flow. We thus compute evolution equations for each coefficient.

We begin with the $xx$-component. Specializing equation \eqref{rf2} to this case, we find
\begin{gather*}
\del_t \{ v_{xx} \}_{1} = \{ \del_t v_{xx}\}_{1}
 =\big\{ 2 (n-1) \big(|{\rm d}x|^2_{\gbar} -1\big) \olg_{xx} \\
 \hphantom{\del_t \{ v_{xx} \}_{1} =}{} -2 ( (n-2) x \hess_{\olg} (x)_{xx} + (x \Delta_{\olg} x) \olg_{xx} ) -2 x^2 \Rc(\gbar)_{xx}\big\}_1 \\
 \hphantom{\del_t \{ v_{xx} \}_{1}}{} =\{ 2 (n-1) \big(|{\rm d}x|^2_{\gbar} -1\big) \olg_{xx} -2 ( (n-2) x \hess_{\olg} (x)_{xx} + (x \Delta_{\olg} x) \olg_{xx} )\}_1,
\end{gather*}
since $x^2 \Rc_{xx}$ already vanishes to second order. Using Lemmas \ref{lemma:expansion} and \ref{appendix:first-odd-term-olgxx} in the appendix we find
\begin{gather*}
\big\{ 2 (n-1) \big(|{\rm d}x|^2_{\gbar} -1\big) \olg_{xx}\big\}_{1} = - 2(n-1) \{v_{xx}\}_1.
\end{gather*}
Using equations \eqref{eqn:hess-asymp-xx} and \eqref{eqn:laplacian-asymp-xx} and Lemma~\ref{lemma:expansion} we obtain
\begin{gather*}
\big\{{-}2 ( (n-2) x \hess_{\olg} (x)_{xx} + (x \Delta_{\olg} x) \olg_{xx} )\big\}_1 = (n-1) \{v_{xx}\}_1 - h_0^{\alpha \beta} \{v_{\alpha \beta}\}_{1}
\end{gather*}
so that
\begin{gather*}
\del_t \{ v_{xx} \}_{1} = -(n-1) \{v_{xx}\}_1 - h_0^{\alpha \beta} \{v_{\alpha \beta}\}_{1}.
\end{gather*}

Passing to the $\alpha \beta$-components we find
\begin{gather*}
\del_t \{ v_{\alpha \beta} \}_{1} =\big\{ 2 (n-1) \big(|{\rm d}x|^2_{\gbar} -1\big) \olg_{{\alpha \beta}} \\
\hphantom{\del_t \{ v_{\alpha \beta} \}_{1} =}{} -2 ( (n-2) x \hess_{\olg} (x)_{{\alpha \beta}} + (x \Delta_{\olg} x) \olg_{{\alpha \beta}} ) -2 x^2 \Rc(\gbar)_{{\alpha \beta}}\big\}_1 \nonumber \\
\hphantom{\del_t \{ v_{\alpha \beta} \}_{1}}{} =\big\{ 2 (n-1)\big(|{\rm d}x|^2_{\gbar} -1\big)\olg_{{\alpha \beta}} -2 ( (n-2) x \hess_{\olg} (x)_{{\alpha \beta}} + (x \Delta_{\olg} x) \olg_{{\alpha \beta}} )\big\}_1.
\end{gather*}

Now using Lemmas \ref{lemma:expansion} and \ref{appendix:first-odd-term-olgxx} once more, we have
\begin{gather*}
\big\{ 2 (n-1) \big(|{\rm d}x|^2_{\gbar} -1\big) \olg_{{\alpha \beta}} \big\}_1 = -2(n-1)\{v_{xx}\}_1 h^0_{\alpha \beta},
\end{gather*}
and using now equations \eqref{eqn:hess-asymp-ab} and \eqref{eqn:laplacian-asymp-ab} we find
\begin{gather*}
\big\{{-}2 ((n-2) x \hess_{\olg} (x)_{{\alpha \beta}} + (x \Delta_{\olg} x) \olg_{{\alpha \beta}} ) \big\}_1 = \{v_{xx}\}_1 - (n-2) \{v_{\alpha \beta}\}_1 + h_0^{\mu \nu} \{v_{\mu\nu}\}_1 h^0_{\alpha \beta}.
\end{gather*}
The evolution equation for $\{v_{\alpha \beta}\}_1$ is therefore
\begin{gather*}
\del_t \{v_{\alpha \beta}\}_1 = (-2n+3)\{v_{xx}\}_1 h^0_{\alpha \beta} - (n-2) \{v_{\alpha \beta}\}_1 - h_0^{\mu \nu} \{v_{\mu\nu}\}_1 h^0_{\alpha \beta}.
\end{gather*}

Thus, the $xx$ and $\alpha\beta$ components satisfy a coupled system of linear equations
\begin{gather*}
\del_t \{ v_{xx} \}_{1} = -(n-1) \{v_{xx}\}_1 - h_0^{\alpha \beta} \{v_{\alpha \beta}\}_{1}, \qquad \mbox{and}\\
\del_t \{v_{\alpha \beta}\}_1 = (-2n+3)\{v_{xx}\}_1 h^0_{\alpha \beta} - (n-2) \{v_{\alpha \beta}\}_1 - h_0^{\mu \nu} \{v_{\mu\nu}\}_1 h^0_{\alpha \beta}.
\end{gather*}
Since each component is zero initially, this condition persists under the flow.

We now specialize to $x\alpha$ indices, and we compute the evolution equation for the coefficient at order $2$, i.e.,
\begin{gather*}
\del_t \{ v_{x\alpha} \}_{2}
 =\big\{ 2 (n-1) \big(|{\rm d}x|^2_{\gbar(t)} -1\big) \olg(t) \\
 \hphantom{\del_t \{ v_{x\alpha} \}_{2}=}{} -2 ( (n-2) \hess_{\olg(t)} (x) + (\Delta_{\olg(t)} x) \olg(t) ) x -2 \Rc(\gbar(t)) x^2 \big\}_2.
\end{gather*}
Now both $\big\{ \olgin_{x\alpha} \big\}_0 = 0$ and $\big\{\Rc(\olgin)_{x\alpha} \big\}_0 = 0$. This forces $\big\{ {-}2 \Rc(\gbar(t)) x^2 \big\}_2 = 0$. Additionally, by the parity property just proved for~$v_{xx}$, i.e., that it remains even to order $2$ we may conclude
\begin{gather*}
\big\{2 (n-1) \big(|{\rm d}x|^2_{\gbar(t)} -1\big) \olg_{x\alpha} \big\}_{2} = 2(n-1)\{\olg^{xx}\}_1 \{ v_{x\alpha}\}_1 = 0
\end{gather*}
as well. This leaves
\begin{gather*}
\del_t \{ v_{x\alpha} \}_{2} = -2 \{ (n-2) x\hess_{\olg(t)} (x)_{x\alpha} + (x \Delta_{\olg(t)} x) \olg(t)_{x\alpha} \}_2.
\end{gather*}

Now, applying equation \eqref{eqn:mixed-ref-1}, we find
\begin{gather*}
-2 \{ (n-2) x\hess_{\olg(t)} (x)_{x\alpha} \}_{2} = (n-2)\partial_{\alpha} \{v_{xx}\}_1 + (n-2)\{ v^{x \mu} \}_1 \{ v_{\mu \alpha}\}_1,
\end{gather*}
and applying equation\eqref{eqn:mixed-ref-2}, we find
\begin{gather*}
\{ -2 (x\Delta_{\olg(t)} x) \olg(t)_{x\alpha} \}_2 = \{\gbar_{x \alpha}\}_1 \big( \{v_{xx}\}_1 - h_0^{\mu \nu} \{v_{\mu\nu}\}_1\big).
\end{gather*}
However, since we already know that $v_{xx}$ and $v_{\alpha \beta}$ are even to order~$2$, the right hand side of both of these equations vanishes. We thus conclude
\begin{gather*}
\del_t \{ v_{x\alpha} \}_{2} = 0.
\end{gather*}
Since this coefficient vanishes initially, $\{v_{x\alpha}\}_2 = 0$ along the flow. We conclude $\olg(t)$ is even to order~$2$ and this concludes the proof of the base case.

\subsection*{Inductive step}

We now assume for the purposes of induction $\olg(t) = \olgin + v$ is even to order $2j$ where $2 \leq 2j \leq 2\ell -2$, and we will prove that~$\olg(t)$ is even to order $2j+2$. By the extended notion of even metric, our inductive hypothesis means that the components
\begin{gather*}
 v_{xx}, \ v_{\alpha \beta}\ \mbox{are even to order} \ 2j, \qquad
 v_{x\alpha} \ \mbox{is odd to order} \ 2j+1,
\end{gather*}
and we must show
\begin{gather*}
v_{xx}, \ v_{\alpha \beta} \ \mbox{are even to order} \ 2j+2, \qquad
v_{x\alpha} \ \mbox{is odd to order} \ 2j+3.
\end{gather*}
The assumption that $\olgin$ is even to order $2\ell$ remains in force. We remark that the inductive hypothesis is used extensively in this proof when extracting the coefficient that breaks parity in an expansion, via Lemma~\ref{lemma:expansion} and the computations of the appendix.

There are essentially no new ideas in the inductive step. Similar to the base case, we first calculate evolution equations for the $(2j+1)$-th coefficients of $v_{xx}$ and $v_{\alpha \beta}$ and again we will find that these coefficients satisfy coupled linear ordinary differential equations in time with zero initial conditions and thus remain zero along the flow. We then use this information and calculation the $(2j+2)$-th coefficient of $v_{x\alpha}$, and find that these coefficients again ordinary differential equations in time with zero initial conditions and thus remain zero along the flow.

Now,
\begin{gather*}
 \del_t \{ v(t) \}_{2j+1} =\big\{ 2 (n-1) \big(|{\rm d}x|^2_{\gbar(t)} -1\big) \olg(t) \\
 \hphantom{\del_t \{ v(t) \}_{2j+1} =}{} -2 ( (n-2) \hess_{\olg(t)} (x) + (\Delta_{\olg(t)} x) \olg(t) ) x -2 \Rc(\gbar(t)) x^2 \big\}_{2j+1}.
\end{gather*}

Specializing first to the $xx$-components, we find first that
\[ \big\{ 2 (n-1) \big(|{\rm d}x|^2_{\gbar(t)} -1\big) \olg(t)\big\}_{2j+1} = -2(n-1) \{v_{xx}\}_{2j+1}. \]
Next, note that
\begin{gather*}
 \big\{ {-}2 \big( (n-2) \hess_{\olg(t)} (x) + (\Delta_{\olg(t)} x) \olg(t) \big) x\big\}_{2j+1} \\
\qquad{} = - 2\left(\!- \frac{1}{2} (n-2)(2j+1) \{ v_{xx} \}_{2j+1} - \frac{1}{2} (2j+1) \{ v_{xx} \}_{2j+1} + \frac{1}{2} (2j+1) h^{\alpha \beta}_0 \{v_{\alpha \beta}\}_{2j+1} \!\right) \\
 \qquad{} = (n-1) (2j+1) \{ v_{xx} \}_{2j+1} - (2j+1) h^{\alpha \beta}_0 \{v_{\alpha \beta}\}_{2j+1}.
\end{gather*}
Finally,
\begin{gather*}
\big\{ {-}2 \Rc(\gbar(t)) x^2 \big\}_{2j+1} = 2j( 2j+1) h^{\alpha \beta}_0 \{v_{\alpha \beta}\}_{2j+1}.
\end{gather*}
We conclude
\begin{gather*}
\del_t \{ v_{xx} \}_{2j+1} = (n-1)(2j-1) \{v_{xx}\}_{2j+1} + (2j-1)(2j+1) h^{\alpha \beta}_0 \{v_{\alpha \beta}\}_{2j+1}.
\end{gather*}

In a similar way, we specialize to the $\alpha \beta$-components
\begin{gather*}
 \del_t \{ v_{\alpha \beta} \}_{2j+1} = -2(n-1) \{v_{xx}\}_{2j+1} h^0_{\alpha \beta} - (n-2)(2j+1) \{ v_{\alpha \beta}\}_{2j+1} \nonumber \\
\hphantom{\del_t \{ v_{\alpha \beta} \}_{2j+1} =}{} + (2j+1) \{v_{xx}\}_{2j+1} h^0_{\alpha \beta} - (2j+1) h_0^{\mu \nu} \{v_{\mu\nu}\}_{2j+1} h^0_{\alpha \beta} + 2j(2j+1) \{v_{\alpha \beta}\}_{2j+1}\\
\hphantom{\del_t \{ v_{\alpha \beta} \}_{2j+1}}{} = (-2(n-1) + (2j+1)) \{v_{xx}\}_{2j+1} h^0_{\alpha \beta} - (2j+1) h_0^{\mu \nu} \{v_{\mu\nu}\}_{2j+1} h^0_{\alpha \beta}\\
\hphantom{\del_t \{ v_{\alpha \beta} \}_{2j+1} =}{} + ( -(n-2) + 2j )(2j+1) \{ v_{\alpha \beta}\}_{2j+1}.
\end{gather*}

These equations form a coupled linear system for the coefficients $\{v_{xx}\}_{2j+1}$ and $\{ v_{\alpha \beta} \}_{2j+1}$, thus these coefficients remain zero along the flow since the initial condition vanishes.

It remains to study the evolution of $\{v_{x \alpha}\}_{2j+2}$, using the improved parity just proved for the~$xx$ and $\alpha \beta$-components. In particular, note that as $\gbar_{x \alpha}$ is odd to order $2j+1$ and $|{\rm d}x|^2_{\gbar}-1$ vanishes to second order and is even to order~$2j+2$,
\[ \big\{ 2 (n-1) \big(|{\rm d}x|^2_{\gbar} -1\big) \olg_{x\alpha}\big\}_{2j+2} = 0. \]

In a similar way, due to equations \eqref{eqn:mixed-ref-1} and \eqref{eqn:mixed-ref-2}, we find that since the $xx$ and $\alpha \beta$-components of $v$ are even to order $2j+2$,
\begin{gather*} \big\{ {-}2 ( (n-2) \hess_{\olg} (x)_{x\alpha} + (\Delta_{\olg} x) \olg_{x \alpha} ) x\big\}_{2j+2} = 0.
\end{gather*}

Finally, equations \eqref{eqn:mixed-ref-3}, \eqref{eqn:mixed-ref-4} and \eqref{eqn:mixed-ref-5} and the parity and coefficients for the Christoffel symbols show that the $2j+2$ component of $-2 x^2 Rc_{x\alpha}$ involves only~$xx$ and $\alpha \beta$ components of~$v$ at order~$2j+1$. These once again vanish. Thus we have proved
\[ \partial_t \{v_{x \alpha}\}_{2j+2} = 0, \]
which completes the induction.

This concludes the proof of Proposition~\ref{prop:even-expansions-preserved}.
\end{proof}

\section{The volume renormalizability condition} \label{VRandRF}

We have now proved that if the initial metric $\gin$ is even to order $2m-2$ (where $2m = \dim M$), then $g(t)$ remains even to order $2m-2$ throughout the interval of existence. Recall that $\gin$ is volume renormalizable when $\tr^{h_0} h_{2m-1} = 0$. In this section we prove
\begin{Proposition} \label{prop:volume-renorm-preserved}If $\gin$ is volume renormalizable, then the solution $g(t)$ to the Ricci flow with initial condition $\gin$ remains volume renormalizable so long as this solution is defined.
\end{Proposition}

The proof proceeds very much as before. We show that some version of the scalar function $\tr^{h_0} h_{2m-1}$ satisfies a homogeneous ordinary differential equation in $t$ with initial condition~$0$, and hence vanishes for $t \geq 0$. The actual quantity we study is slightly more complicated, and is defined below.

Having fixed the initial metric $\gin$ and representative $h_0 \in \frakc\big(\gin\big)$, let $x$ be the corresponding special boundary defining function. We define the $t$-independent metric
\begin{gather*}
G = \frac{{\rm d}x^2 + h_0}{x^2}
\end{gather*}
by truncating the expansion of the tangential metric $h(x)$ in the normal form for~$\gin$. We also set $\olG = x^2 G$.

Next define the function
\begin{gather*}
F := \tr^{\olG} \olg(t) = \tr^{G} g(t).
\end{gather*}
Clearly $F$ is smooth up to $x=0$ for every $t \geq 0$. Moreover, since $G$ is an even metric and $g$ is even to order $2m-2$ by Proposition~\ref{prop:even-expansions-preserved}, it follows $F$ is even to order $2m-2$. Our interest is in the first odd term in the expansion of~$F$.

\begin{Lemma}\label{lemma3.2} Let $\mu = \mu(t,y) := \{F\}_{2m-1}$. Then
\begin{gather*}
\mu = \{\olg_{xx}\}_{2m-1} + (h_0)^{\alpha\beta} \{\olg_{\alpha \beta} \}_{2m-1}.
\end{gather*}
\end{Lemma}
\begin{proof}
Note that since $\olG^{x\alpha} \equiv 0$,
\begin{gather*}
\label{eq3.4}
F = \olG^{xx}\olg_{xx} + \olG^{\alpha \beta} \olg_{\alpha \beta},
\end{gather*}
where all indices are raised with respect to $\olG$. The result is now straightforward.
\end{proof}

\begin{Proposition} \label{prop:linear-system} Setting $\nu = \{\olg_{xx}\}_{2m-1}$ and $\mu = \{\olg_{xx}\}_{2m-1} + (h_0)^{\alpha \beta} \{\olg(t)_{\alpha \beta}\}_{2m-1}$
as above, then
\begin{gather*}
 \partial_t \mu = -2(2m-1) \mu, \qquad
 \partial_t \nu = (2m-3)(2m-1) \mu.
\end{gather*}
\end{Proposition}
\begin{proof}We compute the evolution equation of $F$ using equation \eqref{eq2.1} and substituting $n = 2m$.
\begin{gather*}
\partial_t F = G^{ij} \partial_t g_{ij} = - 2 G^{ij} E_{ij} = -2 \olG^{ij} \overline{E}_{ij} \\
\hphantom{\partial_t F}{} = 2 (2m-1) \big(|{\rm d}x|^2_{\gbar}-1\big) \olG^{ij} \olg_{ij} -2 \olG^{ij} [ (2m-2) x\hess_{\olg} (x)_{ij} +(x\Delta_{\olg} x) \olg_{ij} ]\\
\hphantom{\partial_t F=}{} -2 x^2\olG^{ij}\Rc(\gbar)_{ij} .
\end{gather*}

We must compute the first odd term in this evolution at order $2m-1$, and we split the calculation into three terms.

Since
$F = \olG^{ij} \olg_{ij}$ has leading coefficient $\big(\olG^{ij} \olg_{ij} \big)\vert_{x=0} = 2m$ and is even to order $2m-2$, we have
\begin{gather*}
\big\{ 2 (2m-1) \big(|{\rm d}x|^2_{\gbar}-1\big) \olG^{ij} \olg_{ij}\big\}_{2m-1}
= 2 (2m-1) \big\{\big(|{\rm d}x|^2_{\gbar}-1\big) \olG^{ij} \gbar_{ij} \big\}_{2m-1} \\
\qquad{}= 2 (2m-1) \big( \big\{|{\rm d}x|^2_{\gbar}-1\big\}_0 \big\{ \olG^{ij} \gbar_{ij}\big\}_{2m-1} + \big\{|{\rm d}x|^2_{\gbar}-1\big\}_{2m-1}
\big\{ \olG^{ij} \gbar_{ij}\big\}_{0}\big) \\
\qquad{} = 4m (2m-1) \{\gbar^{xx} \}_{2m-1} = -4m(2m-1) \nu ,
\end{gather*}
where the last equality follows from Lemma \ref{appendix:first-odd-term-olgxx} in the appendix.

Regarding the Hessian terms, we compute applying the calculations given in Appendix~\ref{sec:hess-lap-asymp}:
\begin{gather*}
\big\{{-}2 \olG^{ij}\big[ (2m-2) x\hess_{\olg} (x)_{ij} +(x\Delta_{\olg} x) \olg_{ij} \big]\big\}_{2m-1} \\
\qquad{} = -2 (2m-2) \big\{ \olG^{ij} x\hess_{\olg} (x)_{ij}\big\}_{2m-1} -2 \big\{ (x\Delta_{\olg} x) \olG^{ij} \gbar_{ij}\big\}_{2m-1}\\
\qquad{} = -2 (2m-2) \big\{ x\hess_{\olg} (x)_{xx} + h_0^{\alpha \beta} x\hess_{\olg} (x)_{\alpha \beta}\big\}_{2m-1} -4m \{ x\Delta_{\olg} x\}_{2m-1} \\
\qquad{} = -2( 2m-2) \left( -\frac{1}{2} (2m-1) \{ -\gbar_{xx}+h^{\alpha \beta}_0 \gbar_{\alpha \beta} \}_{2m-1}\right) \\
\qquad \quad{}- 4m \left(-\frac{1}{2} (2m-1) \{ -\gbar_{xx} + h^{\alpha \beta}_0 \gbar_{\alpha \beta} \}_{2m-1}\right) \\
\qquad{}= ( 2m-2) (2m-1) \big\{ {-}\gbar_{xx}+h^{\alpha \beta}_0 \gbar_{\alpha \beta} \big\}_{2m-1}+ 2m (2m-1) \big\{ {-}\gbar_{xx} + h^{\alpha \beta}_0 \gbar_{\alpha \beta} \big\}_{2m-1} \\
\qquad{}= ( 4m-2) (2m-1) \big\{ {-}\gbar_{xx}+h^{\alpha \beta}_0 \gbar_{\alpha \beta} \big\}_{2m-1} = (4m-2)(2m-1) (\mu - 2\nu).
\end{gather*}

Finally, we consider the contributions from the Ricci curvature term. We apply equations~\eqref{eqA.5} and~\eqref{eqA.7}:
\begin{gather*}
\big\{{-}2 x^2\olG^{ij}\Rc(\gbar)_{ij}\big\}_{2m-1} = -2 \big\{ x^2 \Rc(\gbar)_{xx} + x^2 h^{\alpha \beta}_0 \Rc(\gbar)_{\alpha \beta} \big\}_{2m-1} \\
\qquad {} = 2(m-1)(2m-1) h^{\alpha \beta}_0 \{\gbar_{\alpha \beta}\}_{2m-1} + 2(m-1)(2m-1) \big\{ h^{\alpha \beta}_0\gbar_{\alpha \beta}\big\}_{2m-1}\\
\qquad{} = 4(m-1)(2m-1) h^{\alpha \beta}_0 \{\gbar_{\alpha \beta}\}_{2m-1} = 4(m-1)(2m-1) (\mu - \nu).
\end{gather*}

Collecting all of this information, we obtain
\begin{gather*}
\mu' = -2 (2m-1) \mu.
\end{gather*}

To obtain the evolution equation for $\nu = \{\olg_{xx}\}_{2m-1}$, recall that
\begin{gather*}
\partial_t \olg_{xx} = x^2 \partial_t g_{xx} = - 2 \olE( g )_{xx}.
\end{gather*}
Now, equation \eqref{eq2.1} specializes to
\begin{gather*}
-2\olE( g )_{xx} =2(2m-1) \big(|{\rm d}x|^2_{\gbar}-1\big) \olg_{xx} -2 [ (2m-2) x\hess_{\olg} (x)_{xx} +(x\Delta_{\olg} x) \olg_{xx} ]\\
\hphantom{-2\olE( g )_{xx} =}{} -2x^2 \Rc(\gbar)_{xx} .
\end{gather*}
Computing as above gives
\begin{gather*}
\nu' = (2m-1)(2m-3) \mu.\tag*{\qed}
\end{gather*}\renewcommand{\qed}{}
\end{proof}

Applying the result of Proposition \ref{prop:linear-system} we immediately obtain
\begin{Corollary} \label{corollary3.4}
If $(\mu(0), \nu(0)) = (0,0)$ then $(\mu(t), \nu(t)) = (0,0)$ along the flow. In particular, $h_0^{\alpha \beta} \{\gbar_{\alpha \beta}\}_{2m-1} = 0$ for $0 \leq t < T_0$.
\end{Corollary}

We now prove Proposition \ref{prop:volume-renorm-preserved}.

\begin{proof}[Proof of Proposition \ref{prop:volume-renorm-preserved}] Suppose that $\gin$ is a volume renormalizable metric on $M$ with $\dim M = n = 2m$. Let $h_0 \in \mathfrak c\big(\gin\big)$ and $x$ the corresponding special boundary defining function, so that
\begin{gather*}
\gin = \frac{{\rm d}x^2 + h_x}{x^2},
\end{gather*}
where $h_x|_{x=0} = h_0$, and $h_x$ is even to order $2m-2$. As usual, fix tangential coordinates $y^{\alpha}$ extended to be constant along the integral curves of $\nabla^{\olg} x$.

Let $g(t)$ satisfy the Ricci flow equation. We have already shown that it remains even to order $2m-2$. We write $x_t$ for the evolving special boundary defining function corresponding to $g(t)$ and recall that $\mathfrak c(g(t))$ is constant in~$t$, so that
\begin{gather*}
g(t) = \frac{{\rm d}x_t^2 + h_0 + h_2 (t) x_t^2 + \cdots + h_{2m-2}(t) x_t^{2m-2} + h_{2m-1}(t) x_t^{2m-1} + \cdots}{x_t^2} .
\end{gather*}
We wish to prove that $\tr^{h_0} h_{2m-1}(t) = 0$.

Denoting by $x$ the special bdf for the initial metric $\gin$, we have
\begin{gather*}
x_t = {\rm e}^{\omega_t(x,y)} x,
\end{gather*}
where $\omega_t$ vanishes to second order since the conformal infinity is fixed, and $\omega_t$ is even to order $2m-2$ by Lemma \ref{lemma1.7}.
Computing further shows that
\begin{gather*}
{\rm d}x_t = (x\partial_x \omega_t + 1) {\rm e}^{\omega_t} {\rm d}x + (x \partial_{y^{\alpha}} \omega_t) {\rm e}^{\omega_t} {\rm d}y^{\alpha} ,
\end{gather*}
and
\begin{gather*}
{\rm d}x_t^2 = (x\partial_x \omega_t + 1)^2 {\rm e}^{2\omega_t} {\rm d}x^2 + 2(x\partial_x \omega_t + 1) (x \partial_{y^{\alpha}} \omega_t)
{\rm e}^{2\omega_t} {\rm d}y^{\alpha} {\rm d}x\\
\hphantom{{\rm d}x_t^2 =}{} + (x \partial_{y^{\alpha}} \omega_t) (x \partial_{y^{\beta}} \omega_t) {\rm e}^{2\omega_t} {\rm d}y^{\alpha} {\rm d}y^{\beta}.
\end{gather*}
Inserting this into the expression for $g(t)$ yields
\begin{gather*}
g(t) = ( x\partial_x \omega_t + 1)^2 \frac{{\rm d}x^2}{x^2} + 2( x\partial_x \omega_t + 1)( x \partial_{y^{\alpha}} \omega_t) \frac{{\rm d}y^{\alpha} {\rm d}x}{x^2} \\
\hphantom{g(t) =}{} + \big[ {\rm e}^{-2\omega} (h_0)_{\alpha \beta} + (h_2(t))_{\alpha \beta} x^2 + \cdots + (h_{2m-2}(t))_{\alpha \beta} {\rm e}^{(2m-4)\omega} x^{2m-2}\\
\hphantom{g(t) =}{} + (h_{2m-1}(t))_{\alpha \beta} {\rm e}^{(2m-3)\omega} x^{2m-1} + \cdots + (x \partial_{y^{\alpha}} \omega_t) (x \partial_{y^{\beta}} \omega_t) {\rm e}^{2\omega} \big] \frac{{\rm d}y^{\alpha} {\rm d}y^{\beta}}{x^2} .
\end{gather*}
Now observe that by the proof of Lemma~\ref{lemma1.7}, the coefficient of ${\rm d}x^2$ is even to order $2m-2$, and its first odd coefficient is $
\big\{ ( x\partial_x \omega_t + 1)^2 \big\}_{2m-1} = 2 (2m-1) \{\omega\}_{2m-1}$. The coefficient of~${\rm d}x {\rm d}y^{\alpha}$ is odd to order $2m-1$. The entire coefficient of ${\rm d}y^{\alpha} {\rm d}y^{\beta}$ in square brackets of is even to order $2m-2$. Using the fact $\{{\rm e}^{\omega}\}_0 = 1$, the term in square brackets has leading term $h_0$ while the term at order $2m-1$ equals $(h_{2m-1}(t))_{\alpha \beta} - 2 \{\omega\}_{2m-1}(h_0)_{\alpha \beta}$. Applying Corollary~\ref{corollary3.4} with $\nu = 2(2m-1) \{\omega\}_{2m-1}$ and $\mu = 2(2m-1) \{\omega\}_{2m-1} + (h_0)^{\alpha \beta}( (h_{2m-1}(t))_{\alpha \beta} - 2(2m)\{\omega\}_{2m-1}) $ shows that both $\mu$ and $\nu$ must vanish along the flow, and hence we conclude that $\tr^{h_0} h_{2m-1}(t) = 0$ along the flow.
\end{proof}

By combining Propositions \ref{prop:even-expansions-preserved} and \ref{prop:volume-renorm-preserved}, the proof of Theorem \ref{TheoremA} is now complete as well.

\section{Variation of renormalized curvature functionals} \label{VarofCurvFcl}

We now take up the proof of Theorem \ref{TheoremB}. We begin with a quick review of Riesz renormalization, referring to~\cite{Albin} and~\cite{DGH} for more details, and then study the variation of renormalized curvature integrals.

\subsection{Regularized and renormalized integrals}

Fix a product decomposition $[0,\epsilon_0) \times \del M$ of some neighborhood of $\del M$, with projection on the first factor a fixed defining function.
Recall that a function $u$ on $M$ is said to be polyhomogeneous if, in terms of some (and hence any) boundary defining function $x$,
there is an expansion
\begin{gather*}
u \sim \sum_{j} \sum_{\ell=0}^{N_j} u_{j \ell}(y) x^{\gamma_j} (\log x)^\ell,
\end{gather*}
where $\gamma_j$ is a sequence of complex numbers with $\re \gamma_j \to \infty$ and the coefficient functions $u_{j\ell}$ are $\calC^\infty$ on $\del M$.
For simplicity we assume that $\re \gamma_0 \leq \re \gamma_j$ for all $j$. Here the meaning of $\sim$ is that for every $k \in \mathbb{N}$,
\[ u - \sum_{\re \Gamma_j \leq k} \sum_{\ell=0}^{N_j} u_{j \ell}(y) x^{\gamma_j} (\log x)^\ell \in x^k \calA(\Mbar),
\]
where $\calA(\Mbar)$ is the space of conormal functions on $\Mbar$, i.e., $v \in \calA(\Mbar)$ if
$\big|(x\del_x)^j \del_y^\alpha v\big| \leq C_{j,\alpha}$ for all~$j$,~$\alpha$.

Now suppose that $g$ is volume renormalizable and $u$ is polyhomogeneous near $\partial M$. Define
\begin{gather*}
z \mapsto I(z):=\int_M x^z u \,{\rm d} V_g.
\end{gather*}

Expand the volume form ${\rm d}V_g = x^{-n} \sqrt{ \frac{\det h_x}{\det h_0} } \, {\rm d}x {\rm d}V_{h_0} =: x^{-n} J(x,y)\, {\rm d}x {\rm d}V_{h_0}$, where the Jacobian factor $J(x,y)$ has the expansion $\sum\limits_{k \geq 0} J_k x^k$. Writing the expansion of $u$ as above, we see that~$I(z)$ is holomorphic on $\{z \in \CC\colon \re(z)>n-1- \re \gamma_0\}$. Now consider each term
\begin{gather*}
 \int_M x^{z-n+\gamma_j + i}(\log x)^\ell u_{j\ell}(y) J_i(y) \, {\rm d}x {\rm d}y \\
\qquad{} = \int_0^\epsilon \int_{\del M} x^{z-n+\gamma_j + i}(\log x)^\ell u_{j\ell}(y) J_i(y) \, {\rm d}x {\rm d}y \\
\qquad\quad{} + \int_{x \geq \epsilon} \int_{\del M} x^{z-n+\gamma_j + i}(\log x)^\ell u_{j\ell}(y) J_i(y) \, {\rm d}x {\rm d}y,
\end{gather*}
where $J(x,y) \sim \sum J_i(y) x^i$. The second term on the right is entire in $z$, while the first extends meromorphically with a pole of order $\ell+1$ at $z = n - \gamma_j - i - 1$. Thus we conclude $I(z)$ extends to a meromorphic function in the entire complex plane.

\begin{Definition}\label{definition4.1} The \emph{Riesz regularized} integral of $u$ is the finite part of $I(z)$ at $z=0$, i.e.,
\begin{gather*}
\Renint_M u \, {\rm d}V_g:=\FP I(z).
\end{gather*}
\end{Definition}

Note that when $I(z)$ has a simple pole at $z=0$, $\FP I(z) = \lim\limits_{z\to 0} ( I(z)-\Res I(z) )$.

As an application, suppose that $x$ and $\tilde{x} = x {\rm e}^\omega$ are two boundary defining functions, where $\omega \in \calC^\infty(\overline{M})$, then
$\tilde{x}^z - x^z = ({\rm e}^{z\omega}-1) x^z =: A_j(y,z) x^{j + z}$. Suppose also that $u \sim \sum\limits_{j \geq 0} u_j(y) x^j$, i.e., all $\gamma_j$ are nonnegative integers and each $N_j = 0$. A short calculation now shows that if $I(z)$ and $\tilde{I}(z)$ are the regularizations with respect to these two defining functions, then near $z=0$,
\begin{gather*}
\tilde{I}(z) - I(z) = \int_{\partial M} C(y,0) \, {\rm d}V_{h_{0}} +{\mathcal O}(z) ,
\end{gather*}
where $C(y,z) = \sum\limits_{i+j+k=n-1} u_i(y) A_j(y,z) J_k(y)$. The following is an immediate consequence:
\begin{Proposition} \label{proposition4.2}
Suppose that $g$ and $u$ are both even to order $n-2 = 2m-2$. Then $C(y,0)=\frac12 u_{0} A_0 \tr^{h_{0}} h_{n-1}$. If, on the other hand,
the expansion for $u$ starts with the term $x^p$ $(p$~even$)$, then $C(y,0) = \frac12 u_p A_0 \tr^{h_0} h_{n-1-p}$.
\end{Proposition}
\begin{proof}By the parity assumptions, $C(y,0) = u_0(y) A_0(y,0) J_{n-1}(y)$; all other terms vanish. It remains to show that
\begin{gather*}
J_{n-1} = \frac12 \tr^{h_0} h_{n-1}.
\end{gather*}
However, by Proposition~\ref{proposition1.5}, if $g$ is even to order $2m-2$ then $\omega$ is even to order~$2m$, and by simple calculation
$\sqrt{h}/\sqrt{h_{0}}$ is even to order~$2m$ as well. Thus $({\rm e}^{z\omega-1})\sqrt{h}/\sqrt{h_{0}}$ is even to order $2m-2$ and the coefficient of
$x^{2m-1}$ is $\frac{1}{2} \tr^{h_{0}} h_{n-1}$.
\end{proof}

The key example is $u\equiv 1$, in which case the renormalized integral is the renormalized volume. If $g$ is volume renormalizable, then $\tr^{h_{0}}h_{n-1} = 0$,
which shows that the renormalized volume is well-defined with these hypotheses. A very similar argument shows that

\begin{Corollary}\label{corollary4.3}Suppose the conditions of Proposition~{\rm \ref{proposition4.2}} hold and $\tr^{h_{0}} h_{n-1}=0$. Then the renormalized integral of~$u$ is independent of the choice of conformal representative of the boundary metric.
\end{Corollary}

\subsection{Variations of renormalized integrals}
Let $g(t)$ be a family of metrics on a compact manifold $Z$, and consider the Riemannian curvature functional
\begin{gather*}
I_Z(g):=\int_Z u(g) \, {\rm d}V_g,
\end{gather*}
where $u(g)$ is some scalar quantity defined from the curvature tensor and its covariant derivatives. One is often interested in critical points of this action, which are solutions of the Euler--Lagrange equation. The first step in computing these critical points is the variation of the integrand
\begin{gather*}
\left. \frac{\del }{\del t} \right|_{t=0}\int_Z u(t) \, {\rm d}V_{g(t)} = \int_Z \left. \frac{\del }{\del t} \right|_{t=0} \big(u(t) \, {\rm d}V_{g(t)}\big),
\end{gather*}
followed by integration by parts. We are interested in the corresponding calculation on a conformally compact manifold. In particular, we shall only take the action of partially even volume renormalizable metrics and consider variations amongst such metrics. Thus suppose that~$g(t)$ is volume renormalizable for each $t$ and $\gin=:g$; suppose also that the conformal infinity of $g(t)$ is independent of~$t$, with fixed representative~$h_0$. Let $x_t$ be the corresponding family of special boundary defining functions, and write $x_t = {\rm e}^{\omega_t} x$. Let $u(t)=u(g(t))$ be the scalar quantity associated to~$g(t)$. Set
\begin{gather*}
{\cal L}(g) = \left(\frac{\partial}{\partial t} \bigg \vert_{t=0} u(t) + \frac{u(0)}{2} \tr^g \frac{\partial}{\partial t} \bigg \vert_{t=0} g(t)\right).
\end{gather*}
Writing $x_0 = x$ and $u(0) = u$, then we have

\begin{Theorem} \label{thm:renorm-curv-fcl}With the notation above,
\begin{gather*}
\frac{\partial}{\partial t} \bigg \vert_{t=0} \Renint u(t) \, {\rm d}V_{g(t)} =
\frac{\partial}{\partial t} \bigg \vert_{t=0}\FP \; \int_M x_t^z u(t) \, {\rm d}V_{g(t)} = \FP \int_M x^z {\cal L}(g)\, {\rm d}V_{g} ,
\end{gather*}
where the local expression for $\calL$ is exactly the same as in the compact case.
\end{Theorem}

\begin{proof}We begin by computing
\begin{gather*}
\frac{\partial}{\partial t} \bigg \vert_{t=0}\FP \int_M x_t^z u(t) \, {\rm d}V_{g(t)} = \FP \int_M z x^{z-1}_0 \dot{x_t} u(0) \, {\rm d}V_{\gin} +\FP \int_M x^z_0 \frac{\partial}{\partial t}\bigg \vert_{t=0}\big( u(t) \, {\rm d}V_{g(t)}\big ) \\
\hphantom{\frac{\partial}{\partial t} \bigg \vert_{t=0}\FP \int_M x_t^z u(t) \, {\rm d}V_{g(t)}}{}
= \FP \left ( z\int_M x^{z} \dot{\omega_t} u \, {\rm d}V_{g}\right ) + \FP \int_M x^z \frac{\partial}{\partial t}\bigg \vert_{t=0}\big( u(t) \, {\rm d}V_{g(t)}\big )\\
\hphantom{\frac{\partial}{\partial t} \bigg \vert_{t=0}\FP \int_M x_t^z u(t) \, {\rm d}V_{g(t)}}{}
= \Res \int_M x^{z} \dot{\omega_t} u \, {\rm d}V_{g} + \FP \int_M x^z {\cal L}(g) \, {\rm d}V_{g} ,
\end{gather*}
where $\dot{x_t}:=\frac{\partial}{\partial t}\big \vert_{t=0}x_t$, and similarly for $\dot{\omega}_t$.

Since $\omega$ is even to order $2m$, so is ${\dot \omega}$. In addition, ${\dot \omega}\in{\cal O}\big(x^2\big)$ since the conformal infinity is fixed.
This produces a shift by $2$ in the terms in the expansion of ${\dot \omega}u\, {\rm d}V_g$. This means that poles arising from integrating this expression
appear at odd integers. In particular, there is no pole at $z=0$, which means that
\begin{gather*}
\frac{\partial}{\partial t} \bigg \vert_{t=0}\FP \int_M x_t^z u(t)\,{\rm d}V_{g(t)} = \FP \int_M x^z {\cal L}(g) \,{\rm d} V_{g} .\tag*{\qed}
\end{gather*}\renewcommand{\qed}{}
\end{proof}

\begin{Theorem}[Theorem~\ref{TheoremB}]\label{theorem4.5}
Suppose that $\big(M^n,\gin\big)$, $n = 2m$, is volume renormalizable, and let $g(t)$ be a solution of~\eqref{normRF} with $g(0) = \gin$. Then along the flow,
\begin{gather*}
\frac{{\rm d}}{{\rm d}t} \RenV = -\Renint_M S(g(t)) + n(n-1) \, {\rm d} V_{g(t)}.
\end{gather*}
\end{Theorem}
\begin{proof}Proposition \ref{prop:volume-renorm-preserved} asserts that $g(t)$ remains volume renormalizable along the normalized Ricci flow. Consequently the renormalized volume is defined along the flow.

Recall that on a compact manifold the variation of the volume form $g \mapsto {\rm d}V_g$ equals $\frac{1}{2} \tr^{g(t)} \dot{g}(t)$, and for the normalized Ricci flow,
\begin{gather*}
\frac{1}{2} \tr^{g(t)}\partial_t g = - \tr^{g(t)} E_{ij} = -\left( S(g(t)) + n(n-1) \right).
\end{gather*}

Finally, by Theorem \ref{thm:renorm-curv-fcl},
\begin{gather*}
\frac{{\rm d}}{{\rm d}t} \RenV(M, g(t)) = -\Renint_M S(g(t)) + n(n-1) \,{\rm d}V_{g(t)},
\end{gather*}
as asserted.
\end{proof}

\appendix

\section{Parity computations} \label{appendix}

In this appendix we record a number of parity computations used throughout the paper. Recall that we use the notation that if $f(x,y)$ is smooth and has a series expansion of the form $f(x,y) \sim \sum_j a_j(y) x^j$, then we denote by $\{f\}_{n}$ the $n^{\mathrm{th}}$ coefficient function $a_n(y)$. Note that with respect to $0$-derivatives (i.e., derivatives with respect to $x \partial_x$ and $x \partial_{y^{\alpha}} = x \partial_{\alpha}$) one may check $\{ x \partial_x f\}_n = n \{f \}_{n},$ and $\{ x \partial_{{\alpha}} f \}_n = \partial_{{\alpha}} \{f\}_{n-1}.$

It is also easy to check the following
\begin{Lemma} \label{lemma:expansion} If $f$ and $g$ are smooth functions which are even as functions of $x$ to order $2j$, then the product $fg$ is even to order $2j$, with $\{ fg\}_0 = \{f\}_0\{g\}_0$ and $\{ fg\}_{2j+1} = \{ f \}_{2j+1} \{g\}_0 + \{f\}_0 \{g\}_{2j+1}$.
\end{Lemma}

Referring to the definition of even metric given on page~\pageref{defn-even} we have
\begin{Lemma} \label{appendix:first-odd-term-olgxx} Let $g$ be a metric which is even to order~$2j$. Suppose that $x$ is a defining function for which
$|{\rm d}x|^2_{\olg} = \olg^{xx} \sim 1 + \sum\limits_{k=1}^{j} A_{2k} x^{2k} + A_{2j+1} x^{2j+1}$; then completing $x$ to a coordinate system $(x,y)$, we have that
\begin{gather*}
\olg_{xx} = 1 + \sum_{k=1}^{j} a_{2k} x^{2k} + a_{2j+1} x^{2j+1} ,
\end{gather*}
where $a_{2j+1} = -A_{2j+1}$, i.e., $\{ \gbar_{xx} \}_{2j+1} = - \{\gbar^{xx}\}_{2j+1}$.
\end{Lemma}
\begin{proof}We use the formula
\begin{gather*}
\olg^{xx} \olg_{xx} + \olg^{x\alpha} \olg_{x\alpha} = 1,
\end{gather*}
recalling that since $\olg$ is even to order $2j$, the coefficients $\olg_{\alpha x}$ and $\olg^{\alpha x}$ are odd to order $2j+1$, hence may be ignored since their product do not contribute to an odd term before order~$2j+1$. Writing $\olg_{xx} \sim 1 + \sum\limits_{k=1}^{j} a_{2k} x^{2k} + a_{2j+1} x^{2j+1}$, we see that the coefficient at order~$2j+1$ of the product $\olg^{xx} \olg_{xx}$ equals $a_{2j+1} + A_{2j+1}$. However, this coefficient vanishes, which shows that the term in $\gbar_{xx} $ at order $2j+1$ is $-A_{2j+1}$.
\end{proof}

Now suppose $\olg$ is a metric that is even to order $2j$, i.e., $\olg_{xx}$ and $\olg_{\alpha \beta}$ are even to order $2j$ and $\olg_{x\alpha}$ is odd to order $2j+1$. We record the both the parity and the parity breaking term in the expansions of $0$-derivatives of the various metric components. Note that because of the prefactor $x$ in the derivative, all of these terms vanish at $x=0$; this fact is used frequently in what follows along with Lemma~\ref{lemma:expansion}.
$$
\begin{tabu}{|l|l|l|l|}\hline
\mbox{component} &\mbox{parity to order} &\mbox{coefficient after parity broken} \\
\hline
x \partial_x \olg_{xx} & \mbox{even to order} \ 2j & \{ x \partial_x \olg_{xx} \}_{2j+1} = (2j+1) \{\olg_{xx}\}_{2j+1} \\
x \partial_{\alpha} \olg_{xx} & \mbox{odd to order} \ 2j+1 & \{ x \partial_{\alpha} \olg_{xx} \}_{2j+2} = \partial_{\alpha} \{ \olg_{xx}\}_{2j+1} \\
x \partial_x \olg_{x\mu} & \mbox{odd to order} \ 2j+1 & \{ x \partial_x \olg_{x\mu} \}_{2j+2} = (2j+2) \{\olg_{x\mu}\}_{2j+2} \\
 x \partial_{\nu} \olg_{x\mu} & \mbox{even to order} \ 2j+2 & \{ x \partial_{\nu} \olg_{x\mu} \}_{2j+3} = \partial_{\nu} \{\olg_{x\mu}\}_{2j+2} \\
 x \partial_x \olg_{\alpha \beta} & \mbox{even to order} \ 2j & \{ x \partial_x \olg_{\alpha \beta} \}_{2j+1} = (2j+1) \{ \olg_{\alpha \beta}\}_{2j+1}\\
x \partial_{\nu} \olg_{\alpha \beta}& \mbox{odd to order} \ 2j+1& \{ x \partial_{\nu} \olg_{\alpha\beta} \}_{2j+2} = \partial_{\nu} \{\olg_{\alpha\beta}\}_{2j+1} \\
 \hline
\end{tabu}
$$

\subsection*{Christoffel symbols} We now document the expansion of the Christoffel symbols. Once again each of these vanishes at $x=0$, and it emerges from the computation that any such symbol with an even
number of $x$ components is odd to order $2j+1$. Any symbol with an odd number of $x$ components is even to order $2j$ and \textit{the term breaking parity involves only $xx$ and $\alpha \beta$ components of the metric}. Further, let $\hat\Lambda^{\cdot}_{\cdot \cdot}$ denotes the Christoffel symbols involving only tangential derivatives of tangential components of the metric. The calculations leading to the table below are straightforward but tedious, and any entry with an asterisk is not explicitly needed in the sequel.
$$
\begin{tabu}{|l|l|l|l|}\hline
\mbox{component} &\mbox{parity to order} &\mbox{coefficient after parity broken} \\
\hline
x \olGamma_{xx}^x & \mbox{even to order} \ 2j & \big\{x \olGamma_{xx}^x\big\}_{2j+1} = \frac{1}{2} (2j+1) \{\olg_{xx}\}_{2j+1} \tsep{3pt}\\
 x \olGamma^{x}_{\alpha \beta} & \mbox{even to order} \ 2j & \big\{x \olGamma^{x}_{\alpha \beta} \big\}_{2j+1} = -\frac{1}{2} (2j+1) \{ \olg_{\alpha \beta}\}_{2j+1} \\
x \olGamma_{\alpha x}^x & \mbox{odd to order} \ 2j+1 & \big\{ x \olGamma_{\alpha x}^x \big\}_{2j+2} = \big\{\frac{1}{2}\olg^{xx} x \partial_{\alpha} \olg_{xx}\big\}_{2j+2} + \big\{ \frac{1}{2} \olg^{x\mu} x \partial_x \olg_{\alpha \mu} \big\}_{2j+2} \\
x \olGamma_{\alpha x}^{\gamma} & \mbox{even to order} \ 2j & \big\{x \olGamma_{\alpha x}^{\gamma}\big\}_{2j+1} = \frac{1}{2}(2j+1) h_0^{\gamma \beta} \{ \olg_{\alpha \beta}\}_{2j+1} \\
x \olGamma_{x x}^{\gamma} & \mbox{odd to order} \ 2j+1 & *\\
x \olGamma_{\alpha \beta}^{\gamma} & \mbox{odd to order} \ 2j+1& \big\{x \olGamma_{\alpha \beta}^{\gamma}\big\}_{2j+2} = -\frac{1}{2} \{ \gbar^{\gamma x} \}_1 \{x \partial_x \gbar_{\alpha \beta} \}_{2j+1} + \big\{ x \hat{\Lambda}_{\alpha \beta}^{\gamma}\big\}_{2j+2} \\
 \hline
\end{tabu}
$$

\subsection*{Hessian and Laplacian of $\boldsymbol{x}$}\label{sec:hess-lap-asymp}

Recall that
\begin{gather*}
 [\hess_{\olg} x]_{ij} = \partial_i \partial_j x - \olGamma_{ij}^k \partial_k x. \end{gather*}
Thus, from the previous subsection, since $x$ is a coordinate, both
\begin{gather*} [x \hess_{\olg}(x)]_{xx} = - x \olGamma_{xx}^x \qquad \mbox{and} \qquad [x \hess_{\olg}(x)]_{\alpha \beta} = - x \olGamma_{\alpha \beta}^x
\end{gather*}
are even to order $2j$, and
\begin{gather} \label{eqn:hess-asymp-xx}
\{x \hess_{\olg}(x)_{xx}\}_{2j+1} = -\frac{1}{2} (2j+1) \{\olg_{xx}\}_{2j+1}, \\
\{x \hess_{\olg}(x)_{\alpha \beta}\}_{2j+1} = \frac{1}{2} (2j+1) \{\olg_{\alpha \beta}\}_{2j+1}.\label{eqn:hess-asymp-ab}
\end{gather}
Further,
\begin{gather*}
 [x \hess_{\olg}(x)]_{x\alpha} = -x \olGamma_{x\alpha}^x
\end{gather*}
is odd to order $2j+1$, and
\begin{align}
 \{x \hess_{\olg}(x)_{x\alpha}\}_{2j+2} &= -\left\{\frac{1}{2}\olg^{xx} x \partial_{\alpha} \olg_{xx}\right\}_{2j+2} - \left\{ \frac{1}{2} \olg^{x\mu} x \partial_x \olg_{\alpha \mu} \right\}_{2j+2} \nonumber \\
 &= -\frac{1}{2} \{ x \partial_{\alpha} \olg_{xx} \}_{2j+2} - \frac{1}{2} \{ \olg^{x \mu} \}_1 \{ x \partial_x \olg_{\alpha \mu} \}_{2j+1} \nonumber \\
 &= -\frac{1}{2} \partial_{\alpha} \{\olg_{xx}\}_{2j+1} - \frac{1}{2}(2j+1)\{ \olg^{x \mu} \}_1 \{ \olg_{\alpha \mu}\}_{2j+1},\label{eqn:mixed-ref-1}
\end{align}
where we have used the fact that $x \partial_{\alpha} \olg_{xx}$ vanishes to third order in the second equation above.

We also compute leading components of $x \Delta_{\olg} x \, \gbar_{ij}$. First recall that
\begin{gather*}
x \Delta_{\olg} x = x \gbar^{pq} \hess_{pq} x =-x \gbar^{pq} \olGamma_{pq}^x = -x\gbar^{xx} \olGamma_{xx}^x - x\gbar^{x\alpha} \olGamma_{x\alpha}^x - x\gbar^{\alpha \beta} \olGamma_{\alpha \beta}^x.
\end{gather*}
This expression is even to order $2j+1$, and the middle term is even to order $2j+2$. Thus, using Lemma~\ref{lemma:expansion} and the fact that the product of $x$ with any Christoffel symbol vanishes at $x=0$ we have
\begin{align*}
 \{x \Delta_{\olg} x\}_{2j+1} &= -\{ \gbar^{xx}\}_0 \{ x \olGamma_{xx}^x\}_{2j+1} - \{ \gbar^{\alpha \beta}\}_0 \{ x\olGamma_{\alpha \beta}^x\}_{2j+1} \\
 &= -\frac{1}{2} (2j+1) \{\olg_{xx}\}_{2j+1} +\frac{1}{2} (2j+1) h_0^{\alpha \beta}\{\olg_{\alpha \beta}\}_{2j+1}.
\end{align*}

Assembling the above, we discover
\begin{align} \label{eqn:laplacian-asymp-xx}
\{ x \Delta_{\olg} x \gbar_{xx} \}_{2j+1} &= -\frac{1}{2} (2j+1) \{\olg_{xx}\}_{2j+1} +\frac{1}{2} (2j+1) h_0^{\mu\nu}\{\olg_{\mu\nu}\}_{2j+1}, \\
\{ x \Delta_{\olg} x \gbar_{\alpha \beta} \}_{2j+1} &= \left(-\frac{1}{2} (2j+1) \{\olg_{xx}\}_{2j+1} +\frac{1}{2} (2j+1) h_0^{\mu \nu}\{\olg_{\mu\nu}\}_{2j+1} \right)h^0_{\alpha \beta},\label{eqn:laplacian-asymp-ab}
\end{align}
whereas
\begin{align}
\{ x \Delta_{\olg} x \, \gbar_{x \alpha} \}_{2j+2} &= \{ \olg_{x \alpha}\}_1 \{ x \Delta_{\olg} x\}_{2j+1} \nonumber\\
&= \{ \olg_{x \alpha}\}_1 \left( -\frac{1}{2} (2j+1) \{\olg_{xx}\}_{2j+1} +\frac{1}{2} (2j+1) h_0^{\mu \nu}\{\olg_{\mu \nu}\}_{2j+1}\right). \label{eqn:mixed-ref-2}
\end{align}

\subsection*{The Ricci curvature}

We now perform a similar analysis on the components of the Ricci curvature. We write this out in a bit more detail. First, we have
\begin{gather*}
x^2 \overline{\Rc}_{xx} = x^2 {\overline{Rm}_{\alpha x x}}^{\alpha} = x^2 \partial_{\alpha} \olGamma_{xx}^{\alpha} - x^2 \partial_{x} \olGamma_{\alpha x}^{\alpha} + x\olGamma_{xx}^{s} x \olGamma_{\alpha s}^{\alpha} - x \olGamma_{\alpha x}^{s} x \olGamma_{x s}^{\alpha} \nonumber\\
\hphantom{x^2 \overline{\Rc}_{xx}}{} = x^2 \partial_{\alpha} \olGamma_{xx}^{\alpha} - x^2 \partial_{x} \olGamma_{\alpha x}^{\alpha} + x\olGamma_{xx}^{x} x\olGamma_{\alpha x}^{\alpha} + x\olGamma_{xx}^{\mu} x\olGamma_{\alpha \mu}^{\alpha} - x\olGamma_{\alpha x}^{x} x\olGamma_{x x}^{\alpha} - x\olGamma_{\alpha x}^{\mu} x\olGamma_{x \mu}^{\alpha}\nonumber\\
\hphantom{x^2 \overline{\Rc}_{xx}}{} = x \partial_{\alpha} \big(x \olGamma_{xx}^{\alpha}\big) - x^2 \partial_{x} \olGamma_{\alpha x}^{\alpha} + x\olGamma_{xx}^{x} x\olGamma_{\alpha x}^{\alpha} + x\olGamma_{xx}^{\mu} x\olGamma_{\alpha \mu}^{\alpha} - x\olGamma_{\alpha x}^{x} x\olGamma_{x x}^{\alpha}
 - x\olGamma_{\alpha x}^{\mu} x\olGamma_{x \mu}^{\alpha}\nonumber\\
\hphantom{x^2 \overline{\Rc}_{xx}}{} = x \partial_{\alpha} \big(x \olGamma_{xx}^{\alpha}\big) - x \partial_{x} (x \olGamma_{\alpha x}^{\alpha} ) + x \olGamma_{\alpha x}^{\alpha} + x\olGamma_{xx}^{x} x\olGamma_{\alpha x}^{\alpha} + x\olGamma_{xx}^{\mu} x\olGamma_{\alpha \mu}^{\alpha}\nonumber\\
\hphantom{x^2 \overline{\Rc}_{xx}=}{} - x\olGamma_{\alpha x}^{x} x\olGamma_{x x}^{\alpha} - x\olGamma_{\alpha x}^{\mu} x\olGamma_{x \mu}^{\alpha}.
\end{gather*}
Comparing the terms in this expression with the tables above, we find $\Rc_{xx}$ is even to order $2j$. Moreover using the parity for the Christoffel symbols already discussed and that each vanishes at $x=0$, the first and the final four terms are even to order $2j+2$ and thus do not contribute to the term at order~$2j+1$. Thus
\begin{gather}
\big\{ x^2 \overline{\Rc}_{xx} \big\}_{2j+1} = \big\{{-} x \partial_{x} \big(x \olGamma_{\alpha x}^{\alpha}\big) + x \olGamma_{\alpha x}^{\alpha} \big\}_{2j+1}
= \frac{1}{2} \big( {-}(2j+1)^2 + 2j+1 \big) h_0^{\alpha \beta} \{ \olg_{\alpha \beta} \}_{2j+1} \nonumber\\
\hphantom{\big\{ x^2 \overline{\Rc}_{xx} \big\}_{2j+1}}{} = -j (2j+1) h_0^{\alpha \beta} \{ \olg_{\alpha \beta} \}_{2j+1}.\label{eqA.5}
\end{gather}

The tangential Ricci components are
\begin{gather*}
x^2 \overline{\Rc}_{\alpha \beta} x^2 {\overline{Rm}_{s \alpha \beta}}^{s}
= x^2 \big( \partial_s \olGamma_{\alpha \beta}^s - \partial_{\alpha} \olGamma_{s \beta}^s + \olGamma_{\alpha \beta}^r \olGamma_{sr}^s - \olGamma_{s \beta}^r \olGamma_{\alpha r}^s \big) \nonumber\\
\hphantom{x^2 \overline{\Rc}_{\alpha \beta} x^2 {\overline{Rm}_{s \alpha \beta}}^{s} }{} = x^2 \big( \partial_x \olGamma_{\alpha \beta}^x + \partial_{\mu} \olGamma_{\alpha \beta}^{\mu} - \partial_{\alpha} \olGamma_{x \beta}^x - \partial_{\alpha} \olGamma_{\mu \beta}^{\mu}
+ \olGamma_{\alpha \beta}^x \olGamma_{xx}^x + \olGamma_{\alpha \beta}^x \olGamma_{\mu x}^{\mu} \nonumber\\
\hphantom{x^2 \overline{\Rc}_{\alpha \beta} x^2 {\overline{Rm}_{s \alpha \beta}}^{s} =}{} + \olGamma_{\alpha \beta}^{\mu} \olGamma_{x\mu}^x + \olGamma_{\alpha \beta}^{\mu} \olGamma_{\lambda \mu}^{\lambda} - \olGamma_{x \beta}^x \olGamma_{\alpha x}^x - \olGamma_{\mu \beta}^x \olGamma_{\alpha x}^{\mu} - \olGamma_{x \beta}^{\mu} \olGamma_{\alpha \mu}^x - \olGamma_{\lambda \beta}^{\mu} \olGamma_{\alpha \mu}^{\lambda} \big).
\end{gather*}
This entire expression is even to order $2j$, and moreover all terms of the form $\big(x \olGamma\big)\big(x \olGamma\big)$ are even to order $2j+2$, as are the derivative terms involving tangential partial derivatives. We thus find that
\begin{gather}
\big\{ x^2 \overline{\Rc}_{\alpha \beta} \big\}_{2j+1} = \big\{ x^2 \partial_x \olGamma_{\alpha \beta}^x \big\}_{2j+1}
=\big\{ x \partial_x \big( x \olGamma_{\alpha \beta}^x\big) - x \olGamma_{\alpha \beta}^x \big\}_{2j+1}\label{eqA.7}\\
\hphantom{\big\{ x^2 \overline{\Rc}_{\alpha \beta} \big\}_{2j+1}}{}= -\frac{1}{2} (2j+1)^2 \{\olg_{\alpha \beta}\}_{2j+1} + \frac{1}{2} (2j+1) \{\olg_{\alpha \beta}\}_{2j+1}
 = - j (2j+1)\{\olg_{\alpha \beta}\}_{2j+1} .\nonumber
\end{gather}

Now we consider the mixed component $x^2 \overline{\Rc}_{x\alpha}$:
\begin{gather*}
x^2 \overline{\Rc}_{x \alpha} = x^2 {\overline{Rm}_{s x \alpha}}^{s}
= x^2 \big( \partial_s \olGamma_{x \alpha}^s - \partial_{x} \olGamma_{s \alpha}^s + \olGamma_{x \alpha}^r \olGamma_{sr}^s - \olGamma_{s \alpha}^r \olGamma_{x r}^s \big) \\
\hphantom{x^2 \overline{\Rc}_{x \alpha}}{}= x^2 \big( \partial_x \olGamma_{x \alpha}^x + \partial_{\sigma} \olGamma_{x \alpha}^{\sigma} - \partial_{x} \olGamma_{x \alpha}^x - \partial_{x} \olGamma_{\sigma \alpha}^{\sigma}
 + \olGamma_{x \alpha}^x \olGamma_{xx}^x + \olGamma_{x \alpha}^{\nu} \olGamma_{x\nu}^x + \olGamma_{x \alpha}^x \olGamma_{\sigma x}^{\sigma} + \olGamma_{x \alpha}^{\nu} \olGamma_{\sigma \nu}^{\sigma} \\
\hphantom{x^2 \overline{\Rc}_{x \alpha}=}{} - \olGamma_{x \alpha}^x \olGamma_{x x}^x - \olGamma_{x \alpha}^{\nu} \olGamma_{x \nu}^x
- \olGamma_{\sigma \alpha}^x \olGamma_{x x}^{\sigma}- \olGamma_{\sigma \alpha}^{\nu} \olGamma_{x \nu}^{\sigma} \big) \\
\hphantom{x^2 \overline{\Rc}_{x \alpha}}{}= x^2 \big( \partial_{\sigma} \olGamma_{x \alpha}^{\sigma} - \partial_{x} \olGamma_{\sigma \alpha}^{\sigma}+ \olGamma_{x \alpha}^x \olGamma_{\sigma x}^{\sigma} + \olGamma_{x \alpha}^{\nu} \olGamma_{\sigma \nu}^{\sigma} - \olGamma_{\sigma \alpha}^x \olGamma_{x x}^{\sigma}- \olGamma_{\sigma \alpha}^{\nu} \olGamma_{x \nu}^{\sigma} \big).
\end{gather*}
This coefficient is odd to order $2j+1$, thus we will be interested in the coefficient of the first term which breaks parity, at order $2j+2$. Thus purely from parity considerations we find
\begin{gather}
\big\{ x^2\olGamma_{x \alpha}^x \olGamma_{\sigma x}^{\sigma} + x^2\olGamma_{x \alpha}^{\nu} \olGamma_{\sigma \nu}^{\sigma} - x^2 \olGamma_{\sigma \alpha}^x \olGamma_{x x}^{\sigma}- x^2\olGamma_{\sigma \alpha}^{\nu} \olGamma_{x \nu}^{\sigma} \big\}_{2j+2} \nonumber \\
\qquad{} = \big\{x\olGamma_{x \alpha}^x\big\}_1 \big\{ x\olGamma_{\sigma x}^{\sigma}\big\}_{2j+1}+ \big\{x\olGamma_{x \alpha}^{\nu}\big\}_{2j+1} \big\{ x\olGamma_{\sigma \nu}^{\sigma}\big\}_1 \nonumber \\
\qquad\quad {}- \big\{x \olGamma_{\sigma \alpha}^x\big\}_{2j+1} \big\{ x\olGamma_{x x}^{\sigma}\big\}_1 - \big\{x\olGamma_{\sigma \alpha}^{\nu}\big\}_1 \big\{ x\olGamma_{x \nu}^{\sigma}\big\}_{2j+1}.\label{eqn:mixed-ref-3}
\end{gather}
Observe also
\begin{gather}\label{eqn:mixed-ref-4}
\big\{x^2 \partial_{\sigma} \olGamma_{x \alpha}^{\sigma}\big\}_{2j+2} = \big\{ x \partial_{\sigma} \big( x \olGamma_{x \alpha}^{\sigma} \big) \big\}_{2j+2}= \partial_{\sigma}\big\{ x \olGamma_{x \alpha}^{\sigma} \big\}_{2j+1}.
\end{gather}
Finally,
\begin{gather}
\big\{x^2 \partial_{x} \olGamma_{\sigma \alpha}^{\sigma}\big\}_{2j+2} = \big\{x \partial_{x} \big(x\olGamma_{\sigma \alpha}^{\sigma}\big) - x\olGamma_{\sigma \alpha}^{\sigma} \big\}_{2j+2}
 = (2j+1)\big\{x\olGamma_{\sigma \alpha}^{\sigma} \big\}_{2j+2}\nonumber \\
\qquad{} = \frac{-(2j+1)}{2} \{ \gbar^{\sigma x} \}_1 \{x \partial_x \gbar_{\sigma \alpha} \}_{2j+1} + (2j+1)\big\{ x \hat{\Lambda}_{\sigma \alpha}^{\sigma}\big\}_{2j+2} \nonumber\\
\qquad{} = -\frac{(2j+1)^2}{2} \{ \gbar^{\sigma x} \}_1 \{\gbar_{\sigma \alpha} \}_{2j+1} + \frac{1}{2}
\{\gbar^{\sigma \mu}\}_0 \{ x\partial_{\alpha} \gbar_{\mu \sigma}+ x\partial_{\sigma} \gbar_{\alpha \mu} - x\partial_{\mu} \gbar_{\alpha \sigma}\}_{2j+2}.\label{eqn:mixed-ref-5}
\end{gather}

The reader may now verify from equations \eqref{eqn:mixed-ref-3}, \eqref{eqn:mixed-ref-4}, and \eqref{eqn:mixed-ref-5}, that the coefficient $\big\{ x^2 \Rc_{x\alpha} \big\}_{2j+2}$ involves only components $\{ \olg_{xx} \}_{2j+1}$ and $\{ \olg_{\alpha \beta} \}_{2j+1}$ and their tangential derivatives.

\subsection*{Acknowledgements}
EB and RM are grateful to Robin Graham for discussions related to this work. The work of EB was supported by a Simons Foundation grant (\#426628, E.~Bahuaud). The work of EW was supported by an NSERC Discovery Grant RGPIN 203614. RM was supported by the NSF grants DMS-1105050 and DMS-1608223.

\pdfbookmark[1]{References}{ref}
\LastPageEnding


\begin{thebibliography}{99}
\footnotesize\itemsep=0pt

\bibitem{Albin}
Albin P., Renormalizing curvature integrals on {P}oincar\'{e}--{E}instein
 manifolds, \href{https://doi.org/10.1016/j.aim.2008.12.002}{\textit{Adv. Math.}} \textbf{221} (2009), 140--169,
 \href{https://arxiv.org/abs/math.DG/0504161}{arXiv:math.DG/0504161}.

\bibitem{Ammar}
Ammar M., Polyhomog\'en\'eit\'e des m\'etriques compatibles avec une structure
 de {L}ie \`a l'infini le long du flot de {R}icci, Ph.D.~Thesis, Universit\'e
 du Qu\'ebec \`a Montr\'eal, 2019, \href{https://arxiv.org/abs/1907.03917}{arXiv:1907.03917}.

\bibitem{Bahuaud}
Bahuaud E., Ricci flow of conformally compact metrics, \href{https://doi.org/10.1016/j.anihpc.2011.03.007}{\textit{Ann. Inst.
 H.~Poincar\'{e} Anal. Non Lin\'{e}aire}} \textbf{28} (2011), 813--835,
 \href{https://arxiv.org/abs/1011.2999}{arXiv:1011.2999}.

\bibitem{BMW}
Bahuaud E., Mazzeo R., Woolgar E., Renormalized volume and the evolution of
 {APE}s, \href{https://doi.org/10.1515/geofl-2015-0007}{\textit{Geom. Flows}} \textbf{1} (2015), 126--138, \href{https://arxiv.org/abs/1307.4788}{arXiv:1307.4788}.

\bibitem{BahuaudWoolgar}
Bahuaud E., Woolgar E., Asymptotically hyperbolic normalized {R}icci flow and
 rotational symmetry, \href{https://doi.org/10.4310/CAG.2018.v26.n5.a1}{\textit{Comm. Anal. Geom.}} \textbf{26} (2018),
 1009--1045, \href{https://arxiv.org/abs/1506.06806}{arXiv:1506.06806}.

\bibitem{BW}
Balehowsky T., Woolgar E., The {R}icci flow of asymptotically hyperbolic mass
 and applications, \href{https://doi.org/10.1063/1.4732118}{\textit{J.~Math. Phys.}} \textbf{53} (2012), 072501,
 15~pages, \href{https://arxiv.org/abs/1110.0765}{arXiv:1110.0765}.

\bibitem{Bamler}
Bamler R.H., Stability of symmetric spaces of noncompact type under {R}icci
 flow, \href{https://doi.org/10.1007/s00039-015-0317-8}{\textit{Geom. Funct. Anal.}} \textbf{25} (2015), 342--416,
 \href{https://arxiv.org/abs/1011.4267}{arXiv:1011.4267}.

\bibitem{Besse}
Besse A.L., Einstein manifolds, \textit{Ergebnisse der Mathematik und ihrer
 Grenzgebiete~(3)}, Vol.~10, \href{https://doi.org/10.1007/978-3-540-74311-8}{Springer-Verlag}, Berlin, 1987.

\bibitem{Chen-Zhu}
Chen B.-L., Zhu X.-P., Uniqueness of the {R}icci flow on complete noncompact
 manifolds, \href{https://doi.org/10.4310/jdg/1175266184}{\textit{J.~Differential Geom.}} \textbf{74} (2006), 119--154,
 \href{https://arxiv.org/abs/math.DG/0505447}{arXiv:math.DG/0505447}.

\bibitem{CDLS}
Chru\'{s}ciel P.T., Delay E., Lee J.M., Skinner D.N., Boundary regularity of
 conformally compact {E}instein metrics, \href{https://doi.org/10.4310/jdg/1121540341}{\textit{J.~Differential Geom.}}
 \textbf{69} (2005), 111--136, \href{https://arxiv.org/abs/math.DG/0401386}{arXiv:math.DG/0401386}.

\bibitem{DGH}
Djadli Z., Guillarmou C., Herzlich M., Op\'{e}rateurs g\'{e}om\'{e}triques,
 invariants conformes et vari\'{e}t\'{e}s asymptotiquement hyperboliques,
 \textit{Panoramas et Synth\`eses}, Vol.~26, Soci\'{e}t\'{e} Math\'{e}matique
 de France, Paris, 2008.

\bibitem{FG2}
Fefferman C., Graham C.R., The ambient metric, \textit{Annals of Mathematics
 Studies}, Vol.~178, Princeton University Press, Princeton, NJ, 2012.

\bibitem{Gr}
Graham C.R., Volume and area renormalizations for conformally compact
 {E}instein metrics, \textit{Rend. Circ. Mat. Palermo~(2) Suppl.} (2000),
 31--42, \href{https://arxiv.org/abs/math.DG/9909042}{arXiv:math.DG/9909042}.

\bibitem{GL}
Graham C.R., Lee J.M., Einstein metrics with prescribed conformal infinity on
 the ball, \href{https://doi.org/10.1016/0001-8708(91)90071-E}{\textit{Adv. Math.}} \textbf{87} (1991), 186--225.

\bibitem{Guillarmou}
Guillarmou C., Meromorphic properties of the resolvent on asymptotically
 hyperbolic manifolds, \href{https://doi.org/10.1215/S0012-7094-04-12911-2}{\textit{Duke Math.~J.}} \textbf{129} (2005), 1--37,
 \href{https://arxiv.org/abs/math.SP/0311424}{arXiv:math.SP/0311424}.

\bibitem{MM}
Mazzeo R.R., Melrose R.B., Meromorphic extension of the resolvent on complete
 spaces with asymptotically constant negative curvature, \href{https://doi.org/10.1016/0022-1236(87)90097-8}{\textit{J.~Funct.
 Anal.}} \textbf{75} (1987), 260--310.

\bibitem{QSW}
Qing J., Shi Y., Wu J., Normalized {R}icci flows and conformally compact
 {E}instein metrics, \href{https://doi.org/10.1007/s00526-011-0479-7}{\textit{Calc. Var. Partial Differential Equations}}
 \textbf{46} (2013), 183--211, \href{https://arxiv.org/abs/1106.0372}{arXiv:1106.0372}.

\bibitem{R}
Rochon F., Polyhomog\'{e}n\'{e}it\'{e} des m\'{e}triques asymptotiquement
 hyperboliques complexes le long du flot de {R}icci, \href{https://doi.org/10.1007/s12220-014-9505-2}{\textit{J.~Geom. Anal.}}
 \textbf{25} (2015), 2103--2132, \href{https://arxiv.org/abs/1305.5457}{arXiv:1305.5457}.

\bibitem{SSS}
Schn\"{u}rer O.C., Schulze F., Simon M., Stability of hyperbolic space under
 {R}icci flow, \href{https://doi.org/10.4310/CAG.2011.v19.n5.a8}{\textit{Comm. Anal. Geom.}} \textbf{19} (2011), 1023--1047,
 \href{https://arxiv.org/abs/1003.2107}{arXiv:1003.2107}.

\bibitem{WXShi}
Shi W.-X., Ricci deformation of the metric on complete noncompact {R}iemannian
 manifolds, \href{https://doi.org/10.4310/jdg/1214443595}{\textit{J.~Differential Geom.}} \textbf{30} (1989), 303--394.

\bibitem{V}
Vasy A., Microlocal analysis of asymptotically hyperbolic and {K}err--de
 {S}itter spaces (with an appendix by {S}emyon {D}yatlov), \href{https://doi.org/10.1007/s00222-012-0446-8}{\textit{Invent.
 Math.}} \textbf{194} (2013), 381--513, \href{https://arxiv.org/abs/1012.4391}{arXiv:1012.4391}.

\end{thebibliography}
\end{document}